\documentclass[11pt]{amsart}
\usepackage{amsmath,amssymb,enumerate}
\usepackage{latexsym}
\usepackage{amsfonts}
\usepackage {amsbsy}
\usepackage {amsmath}
\usepackage{amssymb}
\usepackage{enumerate}
\usepackage{ bbm}

\newtheorem{theorem}{Theorem}
\newtheorem{lemma}[theorem]{Lemma}
\newtheorem{corollary}[theorem]{Corollary}
\newtheorem{proposition}[theorem]{Proposition}

\theoremstyle{definition}
\newtheorem{definition}[theorem]{Definition}

\newtheorem{remark}[theorem]{Remark}

\numberwithin{theorem}{section}
\numberwithin{equation}{section}

\theoremstyle{definition}

 \numberwithin{equation}{section} %% Comment out for sequentially-numbered
 \numberwithin{figure}{section} %% Comment out for sequentially-numbered
 \theoremstyle{plain}    
 \theoremstyle{plain}    
 \theoremstyle{remark}
 \theoremstyle{remark}
 \theoremstyle{definition}
\theoremstyle{plain}  
\theoremstyle{plain}

\theoremstyle{definition}
\numberwithin{equation}{section}

\newcommand{\B}{\mathbb{B}}
\newcommand{\N}{\mathbb{N}}
\newcommand{\R}{\mathbb{R}}
\newcommand{\C}{\mathbb{C}}

\makeatletter
\@namedef{subjclassname@2010}{%
  \textup{2010} Mathematics Subject Classification}
\makeatother

\begin{document}



\title[Weighted  Green function]{ Weighted  Green functions \\
for complex Hessian operators}

\author{Hadhami El Aini}

\address{Higher School of Sciences and Technology of Hammam Sousse, MAPSFA (LR 11 ES 35), University of Sousse, 4011 Hammam Sousse, Tunisia}

\email{hadhami.elaini@essths.u-sousse.tn}

\author{Ahmed Zeriahi}

\address{Institut de Mathématiques de Toulouse; UMR 5219, Université de Toulouse; CNRS, UPS, 118 route de Narbonne, F-31062 Toulouse Cedex 9, France}

\email{ahmed.zeriahi@math.univ-toulouse.fr}

\date{\today}

\maketitle

\begin{center}
{\it A tribute to Professor Urban Cegrell}
\end{center}

\begin{abstract}
 Let  $1\leq m\leq n$ be two fixed integers. Let $\Omega \Subset \mathbb C^n$ be a bounded  $m$-hyperconvex domain and   $\mathcal A \subset \Omega \times  ]0,+ \infty[$ a finite set of weighted poles. We define and study properties of the $m$-subharmonic Green function of $\Omega$ with prescribed behavior near the weighted set $\mathcal A$. In particular we prove uniform continuity of the exponential Green function in both variables $(z,\mathcal A)$ in the metric space $\bar \Omega \times \mathcal F$, where  $\mathcal F$ is a suitable  family of  sets of weighted poles in $\Omega \times  ]0,+ \infty[$ endowed with the Hausdorff distance. Moreover we give a precise estimate on its modulus of continuity. Our results generalize and improve previous results concerning the pluricomplex Green function du to P. Lelong.
\end{abstract}

%\tableofcontents

\section{Introduction}
Complex Hessian equations have received increasing attention in recent years as they appear in many geometric problems. They provide important examples of fully non-linear PDE’s of second order on complex manifolds which interpolate between (linear) complex Laplace-Poisson equations and (non linear) complex Monge-Ampère equations (see \cite{BZ20} and the references therin).

The pluricomplex Green function ($m=n$) with one pole have been introduced and studied in different contexts by many authors (see \cite{Lem81} \cite{Kli85}, \cite{Dem87}), and  have played an important role in Complex Analysis. Later the pluricomplex Green function with weighted poles was defined and studied by P. Lelong in \cite{Lel89}).

In this paper, we will introduce and study the Green function with weighted poles for the complex Hessian operators, generalizing the pluricomplex Green function with weighted poles considered in \cite{Lel89}. 

Let $\Omega \Subset \C^n$ be a bounded domain and $1 \leq m \leq n$ be a fixed integer. 
Given  a real function $u \in \mathcal{C}^2 (\Omega)$,  for each integer $1 \leq k \leq n$, we denote by $\sigma_k (u)$ the continuous function defined 	at each point $z \in  \Omega$ as the $k$-th symmetric polynomial of the eigenvalues $\lambda_1 (z) \leq \cdots \leq \lambda_n(z)$ 
of the complex Hessian matrix $  \left(\frac{\partial^2 u }{\partial z_j \partial \bar{z}_k} (z)\right)$ of $u$ i.e. 
$$
\sigma_k (u) (z) := \sum_{1 \leq j_1 < \cdots < j_k \leq n} \lambda_{j_1} (z) \cdots \lambda_{j_k} (z).
$$
Recall the usual notations $d = \partial + \bar{\partial}$ and $d^c := (i\slash 2)  ( \bar{\partial} - \partial)$ so  that
$dd^c = i  \partial  \bar{\partial}$.   A simple computation shows that  
$$
(dd^c u)^k \wedge \beta^{n - k} = \frac{(n-k)! \, k!}{n!}  \, \sigma_k (u)  \,  \beta^n, \, \, 
$$
pointwise in $\Omega$ for $1 \leq k \leq m$, where $\beta := dd^c \vert z\vert^2$ is the usual K\"ahler form on $\C^n$.

We say that a real function $u \in \mathcal{C}^2 (\Omega)$ is $m$-subharmonic in $\Omega$ if for any $1 \leq k \leq m$, we have $\sigma_k (u) \geq 0$ pointwise  in $\Omega$.

Observe that  the function $u$ is $1$-subharmonic  in $\Omega$ ($m= 1$) iff it is  subharmonic in $\Omega$ and $\sigma_1 (u) =  \frac{\partial^2 u }{\partial z\partial \bar{z}} = (1\slash 4) \Delta u$, while $u$ is  $n$-subharmonic  in $\Omega$ ($m = n$) iff  $u$ is  plurisubharmonic in $\Omega$ and $\sigma_n (u)  = \mathrm{det}  \left(\frac{\partial^2 u }{\partial z_j \partial \bar{z}_k} \right)$ pointwise  in $\Omega$.

It was shown by Z. B\l ocki  in \cite{Bl05}, that it is possible to define a general notion of $m$-subharmonic function using the concept of $m$-positive currents. Moreover, identifying positive $(n,n)$-currents with positive Radon measures, it is possible to define the $k$-Hessian measure $(dd^c u)^k \wedge \beta^{n - k}$ when $1 \leq k \leq m$ for any (locally) bounded $m$-subharmonic function $u$ on $\Omega$ (see section 2). 

We denote by $\mathcal{SH}_m (\Omega)$ the set of $m$-subharmonic functions in $\Omega$. Then we have
$$
\mathcal{PSH} (\Omega)  = \mathcal{SH}_n (\Omega)\subset \mathcal{SH}_m (\Omega) \subset \mathcal{SH}_1 (\Omega) = \mathcal{SH} (\Omega).
$$
It is possible to extend the Hessian operator $\sigma_m$ to the following class of singular $m$-subharmonic functions :
$$
\mathcal{SH}^b_m (\Omega) := \{u ; u \in \mathcal{SH}_m (\Omega), \exists E \Subset \Omega, u \in L^{\infty} (\Omega \setminus E)\}.
$$
These are $m$-subharmonic functions in $\Omega$ that are bounded near the boundary. The Hessian operator is well defined and continuous under the convergence of decreasing sequences of functions in $\mathcal{SH}^b_m (\Omega)$  (see \cite{Lu12}).  For the complex Monge-Ampère operator this was observed earlier by  J.-P. Demailly \cite{Dem93} and generalized by U. Cegrell \cite{Ceg04}.

The basic example is the fundamental $m$-subharmonic function in $\C^n$ defined as follows :

\begin{equation}\label{eq:DirPb}
\Phi_m (z) = \Phi_{m,n} (z) := \left\{\begin{array}{lcl} 
  - \vert z \vert^{- 2 s },  &\hbox{if}\  1 \leq m < n,  \, \,  \text{with} \, \, \ s := n\slash m - 1\\
  \log (\vert z\vert \slash R_0) & \hbox{if}\  m = n,
\end{array}\right.
\end{equation}
where $\vert \cdot \vert$ is the euclidean norm on $\C^n$ and $R_0 \geq 1$ is large enough so that $\bar{\Omega} \subset \B(0,R_0\slash 2)$.

The function $\Phi_m $  is a radial $m$-subharmonic and negative function in $\C^n$ when $m < n$. It is  plurisubharmonic and negative in $\Omega$ when $m=n$.  Moreover it satisfies the following complex Hessian equation :

\begin{equation} \label{eq:Fsol}
(dd^c \Phi_{m})^m \wedge \beta^{n-m} = c_{m n}  \, \delta_0 \, \beta^n,
\end{equation}
in the sense of currents on $\C^n$, where $\delta_0$ is the unit Dirac measure at the origin and $c_{m n} > 0$ is a numerical constant.

We consider a weight function  $\nu : \Omega \longrightarrow ]0,+ \infty[$ with a finite support $A \subset \Omega$. We associate to this map its graph which is a finite set of {\it weighted poles} in $\Omega$: 

$$
\mathcal A = \{(a,\nu (a)) \ ; \, a \in A\} \subset \Omega \times ]0,+\infty[\}.
$$

 We introduce the Hausdorff distance on the sets of weighted poles. If $\mathcal A , \mathcal A' \subset \Omega \times \R^+$ are two finite sets, we define
 \begin{equation} \label{eq:Hdistance}
 d_H (\mathcal A, \mathcal A'):= \max \{\delta_H (\mathcal A; \mathcal A'), \delta_H (\mathcal A'; \mathcal A)\},
 \end{equation}
 where  
 \begin{equation*}
 \delta_H (\mathcal A; \mathcal A') := \sup_{w \in \mathcal A} d_1 (w , \mathcal A'),
 \end{equation*}
 and   $d_1 \left((a,\nu), (a',\nu')\right) := \vert a - a'\vert + \vert \nu-\nu'\vert$ for $(a,\nu), (a',\nu') \in \Omega \times \R^+$.
 
 It is well known that $d_H$ is actually a distance on the family of all compact subsets of $\Omega \times  ]0,+\infty[$.
 
  We associate to any finite weigthed set of poles $\mathcal A$, its weighted  function
  
\begin{equation} \label{eq:weightfunction}
\phi_m (z,\mathcal A) :=  \inf_{a\in A}  \nu (a) \,  \Phi_{m} (z - a),
\end{equation}
This is a negative but not $m$-subharmonic function in $\Omega$ in general if $m \geq 2$. We also consider the following function

 \begin{equation} \label{eq:weightfunction}
 \psi_m (z,\mathcal A) := \sum_{a \in A} \nu (a) \Phi_m (z-a).
 \end{equation}
 which is $m$-subharmonic and negative in $\Omega$ and satisfies  $\psi_m (\cdot,\mathcal A) \leq \phi_m (\cdot,\mathcal A)$ in $\Omega$.

%Then the $m$-subharmonic envelope of $\phi_m$ in $\Omega$ is well defined as follows
%$$\widehat{\phi}_m (\cdot,\mathcal A) = \widehat{\phi}_m (\cdot,\mathcal A,\Omega) := \sup \{v \in \mathcal{SH}_m (\Omega) ; v \leq \phi_m (\cdot,\mathcal A)\, \, \text{in} \in\, \Omega\}\cdot$$
 %It is  $m$-subharmonic in $\Omega$ and satisfies the following inequalities : 
 %$$\psi_m (\cdot,\mathcal A) \leq \widehat{\phi}_m (\cdot,\mathcal A)  \leq \phi_m (\cdot,\mathcal A), \, \, \text{in} \, \, \, \Omega.$$

We define the associated  $m$-subharmonic Green function of $\Omega$ with weighted poles in $\mathcal A$ by the formula 

$$
G_m (z,\mathcal A) =  G_{m} (z,\mathcal A,\Omega) := \sup \{u (z) ; u \in \mathcal G_m (\Omega,\mathcal A)\}, \, z \in \Omega,
$$
where 

$$
 \mathcal G_{m} (\Omega,\mathcal A) := \{u \in \mathcal{SH}^-_m (\Omega)\, ; \, \exists C_u > 0, u (z) \leq \phi_m (z,\mathcal A) + C_u \, \, \text{in} \, \, \Omega\}.
$$

Observe that $\psi_m \in \mathcal G_{m} (\Omega,\mathcal A) $ so that $G_m (z,\mathcal A) $ is a well defined negative function in $\Omega$.

To state our main results, we need to recall some definitions.   A domain $\Omega \Subset \C^n$ is  said to be $m$-hyperconvex if it admits a negative   exhaustion function $\rho : \Omega \longrightarrow ]-\infty, 0[$ which is $m$-subharmonic in $\Omega$ and
continuous in $\bar \Omega$.  We will say that the domain is  a $m$-hyperconvex domain of Lipschitz type if moreover $\rho$ can be chosen to be Lipschitz continuous in $\bar \Omega$. This terminology is not standard and different from the condition that the domain has a Lipschitz boundary (see Remark \ref{rem:Lipschitz}). 
 \smallskip

 Fix $\delta_0 > 0$ small enough and $0 < \gamma_0 < \gamma_1$.  We  define the following family of sets of weighted poles
  $$
  \mathcal E (\delta_0,\gamma_0, \gamma_1) := \{ \mathcal A \subset \Omega \times \R^+ \, ; \, \delta (\mathcal A) \geq \delta_0, \inf_ {a \in A} \nu(a) \geq \gamma_0, \sum_{a \in A} \nu(a) \leq \gamma_1\},
  $$  
   where $ \delta_{\mathcal A} := \delta(\mathcal A,\partial \Omega)$ is the weighted distance of $\mathcal A$ to the boundary defined by the formula  (\ref{eq:weithtdistance}) below.
   
   \smallskip
   
 Let us state our main results.

\smallskip
\smallskip

{\bf Theorem A.} {\it  Let  $\Omega \Subset \C^n$ be a  $m$-hyperconvex domain and $\mathcal A \subset \Omega \times ]0,+\infty[$ a finite set of weighted poles. Then the associated Green function $G_{m} (\cdot,\mathcal A,\Omega)$ is   $m$-subharmonic and negative in $\Omega$, and satisfies the following properties:

$(1)$ For any $z \in \Omega$, we have
\begin{eqnarray} \label{eq:global-estimate}
\psi_m (z,\mathcal A) \leq  G_{m} (z,\mathcal A,\Omega) 
&\leq  & {\phi}_m (z,\mathcal A) - \Phi_m ( \delta_{\mathcal A}).
\end{eqnarray}

$(2)$  We have the following boundary behaviour
\begin{equation} \label{eq:boundary-estimate}
\lim_{z \to \partial \Omega} \left(\inf_{\mathcal A \in \mathcal  E (\delta_0,\gamma_0,\gamma_1) } G_{m} (z,\mathcal A,\Omega)\right)  = 0.
\end{equation} 

$(3)$  The function $G_{m} (\cdot,\mathcal A,\Omega) \in \mathcal{SH}_m^b (\Omega)$ and satisfies the Hessian equation
\begin{equation} \label{eq:MongeAmperemass}
(dd^c G_{m} (\cdot,\mathcal A,\Omega))^m \wedge \beta^{n-m} = c_{m n} \sum_{a \in A} \nu(a)^m \delta_{a} \beta^n,
\end{equation}
 in the sense of currents on $\Omega$.

$(4)$ The function  $G_{m} (\cdot,\mathcal A,\Omega)$ is the unique $m$-subharmonic function in $\mathcal G_{m} (\Omega,\mathcal A)$ with boundary values $0$, satisfying the complex Hessian equation (\ref{eq:MongeAmperemass}).}

\smallskip
\smallskip

Our second main result gives a precise estimate on the modulus of continuity of the exponential Green function.

\smallskip

For any set of weighted poles $\mathcal A \subset \Omega \times \R^+$ we define the minimal distance between different points in $A$ as follows
\begin{equation} \label{eq:sigma}
2 \sigma_{\mathcal A} :=  \min \{\vert a - b \vert \, ; \, (a,b) \in A^2, a \neq b\}.
\end{equation}
For fixed $\delta_0 > 0$ small enough and $0 < \gamma_0 < \gamma_1$ and $\sigma_0 > 0$,  we define 
$\mathcal F (\delta_0,\gamma_0, \gamma_1,\sigma_0) $ as the family of  sets $\mathcal A \subset \Omega \times \R^+$ satisfying the following conditions :  
\begin{equation} \label{eq:F}
 \delta(\mathcal A,\partial \Omega) \geq \delta_0, \, \, \, \sigma_{\mathcal A} \geq \sigma_0, \, \, \, \inf_{a \in A} \nu (a) \geq \gamma_0, \, \, \, 
\sum_{a \in A} \nu (a) \leq \gamma_1.
 \end{equation}
 Observe that if $\mathcal A  \in \mathcal F (\delta_0,\gamma_0, \gamma_1,\sigma_0) $, then $\mathcal A$ is finite and $\mathrm{Card} (\mathcal A)  \leq \gamma_1 \slash \gamma_0$.
  
 The set $\mathcal F (\delta_0,\gamma_0, \gamma_1,\sigma_0)$ will be endowed with the Hausdorff distance $d_H$ defined by the formula (\ref{eq:Hdistance}).
 
\smallskip

\smallskip

{\bf Theorem B.} {\it Let  $\Omega \Subset \C^n$ be a  $m$-hyperconvex domain of Lipschitz type.
Then the  following properties hold :

1. If $1 \leq m < n$,  for any $0 < \tau <  1 -  \frac{m}{2n - m}$, there exists  constants $M_m$ and $r_1 > 0$ depending on $(\tau, m, n, \delta_0,\gamma_0, \gamma_1) > $ such that for any  $(z',  \mathcal A')  \in \bar \Omega \times \mathcal F (\delta_0,\gamma_0, \gamma_1, \sigma_0)$ and any $(z,  \mathcal A) \in \bar \Omega \times \mathcal F (\delta_0,\gamma_0, \gamma_1,\sigma_0)$, with  $ \vert z' - z\vert + d_H (\mathcal A', \mathcal A) \leq r \leq r_1,$ we have
\begin{equation}  \label{eq:Hldercontinuite}
\exp G_m (z',\mathcal A', \Omega) - \exp G_m (z,\mathcal A, \Omega) \leq M_m r^{\tau}.
\end{equation}

2. If $m=n$, for any $0< \alpha < 1$, there exists  constants $r_1  > 0$, $M_n  > 0$ depending on $(n, \delta_0, \gamma_0,\gamma_1, \alpha)$ such that for any  $(z',  \mathcal A')  \in \bar \Omega \times \mathcal F (\delta_0,\gamma_0, \gamma_1,\sigma_0)$ and any $(z,  \mathcal A) \in \bar \Omega \times \mathcal F (\delta_0,\gamma_0, \gamma_1,\sigma_0)$, with  $ \vert z' - z\vert + d_H (\mathcal A', \mathcal A) \leq r \leq r_1,$  we have
\begin{equation}  \label{eq:uniformcontinuite}
\exp G_n (z',\mathcal A', \Omega) - \exp G_n (z,\mathcal A, \Omega) \leq \frac{M_n }{(\log R_1\slash r)^\alpha},
\end{equation}
where $R_1 := R_0^{1 \slash \gamma_0}$.

In particular the map $ (z,\mathcal A) \longmapsto \exp G_m (z,\mathcal A,\Omega)$ is uniformly continuous  in $ \bar \Omega \times \mathcal F (\delta_0,\gamma_0, \gamma_1,\sigma_0)$.}

\smallskip
Observe that the uniform continuity on the product space $\bar \Omega \times \mathcal F (\delta_0,\gamma_0, \gamma_1,\sigma_0)$  is understood in the sense of the product distance of the euclidean distance on $\bar \Omega$ and the Hausdorff distance $d_H$ on $\mathcal F (\delta_0,\gamma_0, \gamma_1,\sigma_0)$. 
 Let us emphasize that this result is new even in the case of the pluricomplex Green function $(m = n)$, considered by Pierre Lelong  (see \cite{Lel89}). Indeed even in the case of the pluricomplex Green function with one pole, we give a precise estimate of the modulus of continuity of the exponential of the Green function, while Lelong proved only its uniform continuity.
 Moreover in the case  of several weighted poles, we use the Hausdorff distance between the sets of weighted poles rather than the euclidean distance (see Remark \ref{rem:orderedpoles}).

 \section{Preliminaries}
 
 In this section, we recall the basic properties of $m-$subharmonic functions and some known results we will use  throughout the paper. 
 
\subsection{Hessian potentials}
 For a hermitian $n \times n$ matrix $a = (a_{j,\bar k})$ with complex coefficients, we denote by $\lambda_1, \cdots \lambda_n$ the eigenvalues of the matrix $a$. For any $1 \leq k \leq n$ we define the $k$-th trace of $a$ by the formula

$$
S_k (a) := \sum_{1 \leq j_1 < \cdots < j_k \leq n} \lambda_{j_1} \cdots \lambda_{j_k},
$$
which is the $k^{th}$ elementary symmetric polynomial of the eigenvalues $(\lambda_1, \cdots, \lambda_n)$ of $a$.

 Let $\C^n_{(1,1)} $ be the space of real $(1, 1)$-forms on $\C^n$  with constant
coefficients, and define the cone of $m$-positive  $(1,1)$-forms on $\C^n$ by

$$
\Theta_m := \{\omega \in \C^n_{(1,1)}  \, ; \,  \omega \wedge  \beta^{n - 1} \geq 0, \cdots,  \omega^m \wedge  \beta^{n - m} \geq 0\}.
$$

\begin{definition}
1) A smooth $(1,1)$-form $\omega$ on $\Omega$ is said to be $m$-postive on $\Omega$ if for any $z \in \Omega$, $\omega (z) \in \Theta_m$.

2) A function $u:\Omega \rightarrow \mathbb{R}\cup\{-\infty\}$ is said to be  $m-$subharmonic  on $\Omega$ if it is subharmonic on $\Omega$ (not identically $-\infty$ on any component) and  for any collection of smooth $m-$positive $(1,1)-$forms  $\omega_1,...,\omega_{m-1}$ on $\Omega$, the following inequality holds in the sense of currents 

  $$
  dd^c u\wedge \omega_1\wedge...\wedge \omega_{m-1} \wedge \beta^{n-m}\geq 0,
  $$
  in the sense of currents on $\Omega$.
\end{definition}

We denote by $\mathcal{SH}_m (\Omega) $ the positive convex cone of $m$-subharmonic functions on $\Omega$ which are not identically $-\infty$ on any component of $\Omega$.  These are the $m$-Hessian potentials.

We give below the most basic properties of $m$-subharmonic functions that will be used in the sequel (see \cite{Bl05}, \cite{Lu12}).

\begin{proposition}\label{prop:basic}

\noindent 1.  If $u\in \mathcal{C}^2(\Omega)$, then $u$ is  $m$-subharmonic on $\Omega$ if and only if $(dd^c u)^k\wedge \beta^{n-k}\geq0$
pointwise on $\Omega$ for $k=1, \cdots, m$.

 \noindent 2. $\mathcal{PSH}(\Omega)=\mathcal{SH}_n(\Omega)\subsetneq \mathcal{SH}_{n-1}(\Omega)\subsetneq...\subsetneq \mathcal{SH}_1(\Omega)=\mathcal{SH}(\Omega) $.
 
\noindent 3.  $\mathcal{SH}_m(\Omega) \subset L^1_{loc} (\Omega)$ is a positive convex cone. 
  
\noindent 4.  If $u$ is $m$-subharmonic on $\Omega$ and $f: I \rightarrow\mathbb{R}$ is a  convex, increasing function on some interval containing the image of $u$, then $f \circ u$ is $m$-subharmonic on $\Omega$.

\noindent 5. The limit of a decreasing sequence of  functions in $\mathcal{SH}_m(\Omega)$ is $m$-subharmonic on $\Omega$ when it is not identically $- \infty$ on any component.
 
\noindent 6.  Let $u$ be an $m$-subharmonic function in $\Omega$. Let $v$ be an $m$-subharmonic function in a domain $\Omega'  \subset \C^n$ with $\Omega \cap \Omega' \neq \emptyset$. If $u \geq v$ in $\Omega \cap \partial\Omega'$, then the function
   $$
   w(z):=\left\{\begin{array}{lcl}
\max \{u(z),v(z)\}  &\hbox{ if}\ z \in \Omega \cap \Omega'\\
 u(z)  &\hbox{if}\  z \in\Omega\setminus\Omega'\\
\end{array}\right.
$$
is $m$-subharmonic in $\Omega$.
 \end{proposition}
 \smallskip
 \subsection{The comparison principle}
 
The following result is well known (see \cite{Lu12}, \cite{Lu15}).

\begin{proposition} \label{prop:CP1}
 Assume that $u,v\in \mathcal{SH}_m(\Omega)\cap L^{\infty}(\Omega)$ and for any $\zeta \in \partial \Omega$, $\liminf_{z \rightarrow \zeta }(u(z)- v(z))\geq 0$.  Then 
 $$
 \int_{\{u<v\}}(dd^c v)^m\wedge\beta^{n-m} \leq \int_{\{u<v\}}(dd^c u)^m\wedge\beta^{n-m}.
 $$
 Consequently, if $(dd^cu)^m\wedge\beta^{n-m}\leq(dd^cv)^m\wedge\beta^{n-m}$ weakly on $\Omega$, then $u \geq v$ in $\Omega$.
\end{proposition}

As a consequence we can deduce the following result.
 
 \begin{proposition} \label{prop:CP2}  Assume that $u,v\in \mathcal{SH}^b_m(\Omega)$ and for any $\zeta \in \partial \Omega$, $\lim_{z \rightarrow \zeta }(u(z)- v(z)) = 0$.  Then
 
 1) if $u = v $ near the boundary $\partial \Omega$, we have
 $$
 \int_\Omega (dd^c u)^m\wedge\beta^{n-m} =  \int_\Omega (dd^c v)^m\wedge\beta^{n-m} ,
 $$
 
 2) if $u \leq v$ in $\Omega$,
 $$
 \int_\Omega (dd^c u)^m\wedge\beta^{n-m} \geq  \int_\Omega (dd^c v)^m\wedge\beta^{n-m}.
 $$
 \end{proposition}
   Here is another important tool for comparing $m$-subharmonic functions, called the domination principle.
    \begin{proposition} \label{prop:DP}  Let  $u,v\in \mathcal{SH}^b_m(\Omega)$ such that  $\liminf_{z \rightarrow \zeta }(u(z)- v(z))
    \geq 0$, for any $\zeta \in \partial \Omega$. Assume that $ u \geq v$, almost everywhere in $\Omega$ with respect to the Hessian measure $(dd^c u)^m \wedge \beta^{n-m}$. Then $u \geq v$ everywhere in $\Omega$.
    \end{proposition}
    This result was proved in the case $m=n$ by Bedford and Taylor (see \cite[Corollary 4.5]{BT82}) using the comparison principle Proposition \ref{prop:CP1}.
    The same proof is valid in the general case.
    \smallskip
 \subsection{Comparison of residual masses}
 
 We will need the following comparison Theorem inspired by a result of J.-P. Demailly for the complex Monge-Ampère operator (see\cite{Dem93})
\begin{lemma} \label{lem:residualmass}
Let $u,v\in \mathcal {SH}^b_m(\Omega)$ such  that 

$$
\ell :=\limsup\frac{u(z)}{v(z)}< \infty\ \ \ as  \ z\in\Omega, \ v(z)\rightarrow -\infty.
$$

Then, 
$$\int_{\{v = - \infty\}} (dd^{c}u)^m\wedge\beta^{n-m}\leqslant \ell^m  \int_{\{v = - \infty\}}  (dd^cv)^m\wedge\beta^{n-m}.
$$  
In particular if $l=\lim\frac{u}{v}$ as $\ z\in\Omega,  \ v(z)\rightarrow -\infty $, we have equality.
\end{lemma}
\begin{proof}
The proof is the same as  the one for the complex Monge-Ampère operator (see \cite{Dem93}). For convenience we give it here.
It is sufficient it to prove for $l=1$. We can assume that $u \leq 0$ and $v \leq 0$ on a neighborhood of $\{v = - \infty\}$.
 Fix  $c>0$. By assumption given $\varepsilon >0$,  there exists $b >1$ large enough so that $u  - c \geq (1 + \varepsilon) v (z) =: v_\varepsilon$ on the set $\{v (z) < - b\} \Subset \Omega$.
 
 We consider the $m$-subharmonic  function $w_c :=\max \{u-c,v_\varepsilon \}$ on $\Omega$ which satisfies
 $w_c = u - c$ on the open set $\{v (z) < - b\} \Subset \Omega$. Therefore for any $c > 0$ and $\varepsilon > 0$ there exists $b > 1$ such that 
 %we have for any $b > 1$ large enough,
% $$
 %\int_{\{v < - b\}}(dd^{c}w_c)^{m}\wedge\beta^{n-m}=\int_{\{v < - b\}} (dd^{c}u)^{m}\wedge\beta^{n-m},
 %$$
 %hence
 \begin{equation} \label{eq:currents}
  (dd^{c}w_c)^{m}\wedge\beta^{n-m}= (dd^{c}u)^{m}\wedge\beta^{n-m},
 \end{equation}
 in the sense of currents on  the open set $\{v (z) < - b\}$.
 
  On the other for a fixed $\varepsilon >0$, $w_c$ decreases to $v_\varepsilon$ in $\Omega$ as $c$ increases to $+ \infty$. Hence by continuity of the Hessian operator for decreasing sequences of $m$-subharmonic functions in $\mathcal{SH}^b_m(\Omega)$, it follows  that 
  $$
  (dd^{c}w_c)^{m}\wedge\beta^{n-m} \to (dd^{c} v_\varepsilon)^{m}\wedge\beta^{n-m},
  $$
   in the sense of currents on $\Omega$ as c increases to $+ \infty$ (see \cite{Lu15}).
 
 Now fix a compact set $K \subset \{v = - \infty\}$. Then by (\ref{eq:currents}) and upper semi-continuity we have
 \begin{eqnarray*}
\int_K  (dd^{c}u)^{m}\wedge\beta^{n-m} &=  &\limsup_{c \to + \infty}  \int_K (dd^{c} w_c)^{m}\wedge\beta^{n-m}  \\
&\leq & \int_{K} (dd^{c} v_\varepsilon)^{m}\wedge\beta^{n-m} \\
 &=& (1+\varepsilon)^m \int_{K} (dd^{c} v)^{m}\wedge\beta^{n-m}.
 \end{eqnarray*}
  
 Letting $\varepsilon \to 0$ we obtain the inequality
 $$
 \int_K  (dd^{c}u)^{m}\wedge\beta^{n-m} \leq \int_K  (dd^{c}v)^{m}\wedge\beta^{n-m},
 $$
 for any compact subset $K \subset \{v = - \infty\}$.
 Since by definition the two currents extend as positive Borel measures with locally finite mass on $\Omega$, by interior regularity we obtain the same inequality for the Borel set $\{v = - \infty\}$ i.e.
 $$
 \int_{\{v = - \infty\}}  (dd^{c}u)^{m}\wedge\beta^{n-m} \leq \int_{\{v = - \infty\}}  (dd^{c}v)^{m}\wedge\beta^{n-m},
 $$
 which is the required inequality.
 
 Actually the previous proof gives more information, namely we have the following inequality
 $$
{\bf 1}_{ \{v = - \infty\}}  (dd^{c}u)^{m}\wedge\beta^{n-m} \leq {\bf 1}_{ \{v = - \infty\}}  (dd^{c}v)^{m}\wedge\beta^{n-m},
 $$
 in the sense of Borel measures on $\Omega$. Here ${\bf 1}_{ \{v = - \infty\}} $ is the characteristic function of the Borel set 
 $\{v = - \infty\}$.
\end{proof}

\smallskip

\subsection{The maximal sub-extension}

Let $\Omega \Subset \C^n$ be a bounded domain and $D \subset \Omega$ be an open subset. Let  $h : D \longrightarrow [-\infty, 0]$  be an upper semi-continuous function in $D$ (the obstacle function).  A function $u_0 \in \mathcal{SH}^-_m(\Omega)$ is called a  $m$-subharmonic sub-extension of $h$ to $\Omega$ if $u_0 \leq h$ in $D$. If such a sub-extension exists, we can consider the maximal $m$-subharmonic sub-extension of $h$ to $\Omega$ defined in $\Omega$ as follows :
$$
U = U_{D,\Omega} (h) := \sup \{ u \in \mathcal{SH}^-_m(\Omega), u \leq h \, \text{in} \, \, D\}.
$$
This construction is classical in Potential Theory and has been considered also in different contexts in Pluripotential Theory (see \cite{BT76,BT82}, \cite{CKZ11}, \cite{GLZ19}, \cite{BZ20}).

Here we will need the following result.
\begin{proposition} \label{prop:Subext} Let $h : D \longrightarrow [-\infty, 0]$  be an upper semi-continuous function in $D$  which admits a  negative $m$-subharmonic sub-extension $u_0$ to $\Omega$. Then its maximal  $m$-subharmonic sub-extension $U$ to $\Omega$ is $m$-subharmonic in $\Omega$ and satisfies $u_0 \leq U$ in $\Omega$ and  $U \leq h$ in $D$.

Furthermore if $h \in \mathcal{SH}_m (D)$ and $ u_0 \in \mathcal{SH}_m^b (\Omega)$, then $U \in \mathcal{SH}_m^b (\Omega)$ and  the $m$-Hessian measure of  $U$  is carried by the contact set $ \mathcal Q := \{ z \in D ; U (z)= h(z)\}$ i.e. 
$$
\int_{\{ U< h\}} (dd^c U)^m \wedge \beta^{n-m}  = 0,
$$
where $\{ U< h\} := \{z \in D \, ; \, U(z) < h(z)\} = D \setminus \mathcal Q$ is the non-contact set.
\end{proposition}
In the case $m=n$, the result follows essentially from \cite[Theorem 2.1]{CKZ11}.  For a bounded lower semi-continuous obstacle $h$, this was considered in  \cite{GLZ19}.
\begin{proof} If $h$ is a bounded and continuous  function in $D$, the function $U$ is bounded and $m$-subharmonic in $\Omega$ and the set $\{ U< h\} $ is an (euclidean) open subset of $D$. The result can then be easily proved using the classical method of balayage in each ball $B \Subset \{ U< h\} $ to show that $(dd^c U)^m \wedge \beta^{n-m} = 0$ in $B$ (see the proof of Proposition \ref{prop:Hessian} below).
In the general case, the proof of \cite[Theorem 2.1]{CKZ11}  can be easily adapted to our situation using the fact that if $h$ is $m$-subharmonic in $D$, it is quasi-continuous with respect to the $m$-Hessian capacity (see \cite{Lu12,Lu15}).
\end{proof}
 \section{The weighted Green function}

The goal of this section is to prove Theorem A  in several steps.
 \subsection{Global estimates}

We first define  the weigthed radius function : 
 \begin{equation}\label{eq:weigth-distance}
\theta_m (\delta,\nu) :=   \Phi_m^{-1} (\Phi_m(\delta)\slash \nu) = \left\{\begin{array}{lcl} 
  \nu^{1 \slash 2 s} \delta,  &\hbox{if}\  1 \leq m < n,  \\
R_0 (\delta \slash R_0)^{1 \slash \nu} & \hbox{if}\  m = n,
\end{array}\right.
\end{equation}
where $s := n\slash m - 1> 0$.

Observe that the function $\theta_m$ is increasing in each variable  $(\delta,\nu) \in ]0,R_0[ \times ]0,+ \infty[$.
 Next we define  the weighted distance $\delta_{\mathcal A}$ of $\mathcal A$ to the boundary $\partial \Omega$  as follows:
 \begin{equation}\label{eq:weithtdistance}
 \delta ({\mathcal A}, \partial \Omega)  :=  \inf_{(a,\nu) \in \mathcal A} \theta_m (d (a), \nu^{-1}) = \left\{\begin{array}{lcl} 
  \inf_{(a,\nu) \in \mathcal A} \nu^{-1 \slash 2 s} d (a),  &\hbox{if}\  1 \leq m < n,  \\
 \inf_{(a,\nu) \in \mathcal A}  R_0 \left(d (a)\slash R_0\right)^\nu & \hbox{if}\  m = n,
\end{array}\right.
\end{equation}
Here $ d (a) = d(a,\partial \Omega) := \inf \{\vert a - \zeta \vert ; \zeta \in \partial \Omega\}$ is the euclidean distance of $a$  to the boundary of $\Omega$ and $R_0 = 2  \,  \hbox{diam} (\Omega)$.

Observe that  by definition, we have for any $(a,\nu) \in \mathcal A$,
 \begin{equation}\label{eq: weighted-radius}
0 < \delta < \theta_m (d (a), \nu^{-1}) \Longleftrightarrow  \,  \bar B(a,\theta_m (\delta,\nu)) \subset \Omega.
 \end{equation}
 
For  any fixed $\delta > 0$, we define the sublevel set of the weight function as follows :
 $$
 A_{\delta} := \{z \in \C^n \, ; \, \phi_m(z,\mathcal A) < \Phi_m(\delta) \} =  \bigcup_{(x,\nu) \in   \mathcal  A} B(x,\theta_m(\delta,\nu)).
 $$
Then
 \begin{equation}\label{eq: weight-function}
0 < \delta <  \delta_{\mathcal A}  \Longleftrightarrow \forall (a,\nu) \in \mathcal A, \,  \bar B(a,\theta_m (\delta,\nu)) \subset \Omega \Longleftrightarrow A_\delta \Subset \Omega,
 \end{equation}
 and
  \begin{equation}\label{eq: ineqality}
z \notin   B(a,\theta_m (\delta,\nu))  \Longleftrightarrow  \nu \Phi_m (z-a) - \Phi_m (\delta) \geq 0.
 \end{equation}

Finally, recall that the minimal  distance between distinct points in $A$ is defined by 
  \begin{equation}\label{eq: minimaldistance}
  2 \sigma_{\mathcal A} := \min \{\vert a - b\vert \;,  \, (a,b) \in A^2, a \neq b\}.
  \end{equation} 
  
  The following lemma will be crucial.
  \begin{lemma} Fix $0 < \gamma_0 \leq   1  \leq   \gamma_1$ and let $\mathcal A \subset \C^n \times \R^+$ be a finite set such that
  $$
   \inf \{\nu (a), a \in A\} \geq \gamma_0, \, \, \, \, \sum_{a \in A} \nu(a) \leq \gamma_1.
  $$ 
 
  Then the following estimates hold :
  
  1) for any $\delta  >0$ and $z \notin A_{\delta}$, we have
   \begin{equation} \label{eq:FIneq1} 
   \phi_{m}(z,\mathcal A) +  \gamma_1  \gamma_0^{-1} \Phi_m (\delta) \leq \psi_{m}(z, \mathcal A)\leq \phi_m(z,\mathcal A),
\end{equation} 

  2) for any $0< \delta \leq \theta_m(\sigma_{\mathcal A},\gamma_1^{-1})$ and $z \in A_{\delta}$, we have
  \begin{equation} \label{eq:FIneq2} 
   \phi_{m}(z,\mathcal A)  +  \Phi_m (\delta) \leq \psi_{m}(z, \mathcal A)\leq \phi_m(z,\mathcal A),
\end{equation} 
In particular  for any $z \in \C^n$, we have
 \begin{equation} \label{eq:FIneq3}
 \psi_{m}(z, \mathcal A)\leq   \phi_{m}(z, \mathcal A) \leq \psi_{m}(z, \mathcal A)  -   {\gamma}_1^2 \gamma_0^{-1} \Phi_m(\sigma_{\mathcal A}).
\end{equation} 

  \end{lemma}
  
  \begin{proof}
 For conveniency we use the following notation for the sets of weighted $\mathcal A = \{(a,\nu(a) \, ; \,  a \in A\}$. Fix $z \in \Omega$. Then there exists $(a,\nu (a)) \in \mathcal A$ such that for any $(b,\nu(b)) \in \mathcal A$,
  $$
  \phi_m(z,\mathcal A) =  \nu(a) \Phi_m(z-a)  \leq \nu(b) \Phi_m(z-b),
  $$
  and then 
  \begin{equation} \label{eq:Ineq0}
  \psi_m(z,\mathcal A) = \phi_m(z,\mathcal A)  + \sum_{b \neq a}  \nu(b) \Phi_m(z-b)\cdot
  \end{equation}
  Assume first that $z \notin A_{\delta} := \bigcup_{x \in  A} B(x,\theta_m(\delta,\nu (x)))$. Then for any $b \in A$, $\vert z-b\vert \geq \theta_m(\delta,\nu (b)),$ hence  $\nu (b)  \Phi_m(z-b) \geq \Phi_m (\delta)$ and then by (\ref{eq:Ineq0}), we have
  $$
  \psi_m(z,\mathcal A) \geq \phi_m(z,\mathcal A)  + (p -1)  \Phi_m (\delta),
  $$
  where $p$ is the cardinality of $\mathcal A$. Since $p \leq  \gamma_1  \gamma_0^{-1}$, we obtain the  inequality (\ref{eq:FIneq1}).

  Now assume that $z \in A_{\delta} $ then there exists $x \in A$ such that $z \in B(x,  \theta_m(\delta,\nu(x))$, hence  
$\vert z - x\vert \leq    \theta_m(\delta,\nu (x)) \leq  \theta_m(\delta,\gamma_1)$.
  Moreover if  $0 < \delta \leq \theta_m(\sigma_{\mathcal A},\gamma_1^{-1})$, we infer that  $\sigma_{\mathcal A} \geq \theta_m(\delta,\gamma_1)$ and then  for any $b \in A \setminus \{x\}$, we have 
   \begin{eqnarray*} \label{eq:Ineq2}
  \vert z - b\vert \geq  \vert b-x\vert - \vert z-x\vert) &\geq & 2 \sigma_{\mathcal A}  -  \theta_m(\delta,\gamma_1) \nonumber \\
  & \geq & \theta_m(\delta,\gamma_1).
  \end{eqnarray*}
  Hence $ \Phi_m(z-b) \geq  \Phi_m \circ \theta_m (\delta,\gamma_1)$ and then from (\ref{eq:Ineq0}), it follows  that 
   $$
  \psi_m(z,\mathcal A) = \phi_m(z,\mathcal A)  + \sum_{b \neq a}  \nu(b) \Phi_m(z-b) \geq \phi_m(z,\mathcal A)  + \Phi_m(\delta).
  $$ 
  This proves (\ref{eq:FIneq2}) and (\ref{eq:FIneq3}) follows immediately.

 % Now assume that $z \in \B(a, \theta(\nu_a,\delta)$  then 
 % \begin{equation} \label{eq:Ineq1}
 % \nu_a \Phi_m(z-a) \leq \Phi_m(\delta).
 % \end{equation}
 % On the other hand fix  $0 < \delta \leq \gamma_1^{-2s} \sigma_{\mathcal A}$  and  observe that for any $b \in A \setminus \{a\}$, we have 
   %\begin{eqnarray} \label{eq:Ineq2}
  %\vert z - b\vert \geq  \vert b-a\vert - \vert z-a\vert) &\geq & 2 \gamma_1^{2s} \delta  - \nu(a)^{1\slash 2s}\delta \\
  %& \geq & \nu(b)^{1\slash 2s} \delta \geq \gamma_0^{1\slash 2s} \delta.
  %\end{eqnarray}
   %It follows from (\ref{eq:Ineq1}) and (\ref{eq:Ineq2}) that for  $0 < \delta \leq \gamma_1^{-2s} \sigma_{\mathcal A}$  and $z \in \B(a, \nu(a)^{1\slash 2s}\delta)$, we have
   %$$
  %\phi_m(z,\mathcal A) =\frac{- \nu(a)}{\vert z-a \vert^{2s}},
  %$$
   %and 
   %$$
  %\psi_m(z,\mathcal A) = \phi_m(z,\mathcal A)  + \sum_{b \neq a} \frac{- \nu_b}{\vert z-b \vert^{2s}} \geq \phi_m(z,\mathcal A) - \gamma_1 \gamma_0^{-1} \delta^{-2s}.
  %$$ 
  %This proves (\ref{eq:FIneq2}).
  \end{proof}
   As a consequence we have the following useful estimates for the Green function.
  \begin{corollary}  \label{Cor:Fest3} Let $\mathcal A \subset \Omega \times \R^+$ be a finite set. Then  for $z \in \Omega$, we have
  \begin{equation} \label{eq:FIneq4}
\psi_{m}(z, \mathcal A)\leqslant G_m (z,\mathcal A,\Omega)\leqslant \phi_{m}(z,\mathcal A) - \Phi_m(\delta_{\mathcal A}) ,
\end{equation} 
and
  \begin{equation} \label{eq:FIneq5}
 \psi_m(z,\mathcal A) \leq  G_m(z,\mathcal A,\Omega) \leq  \psi_m (z,\mathcal A) - \Phi_m(\delta_{\mathcal A})  -   \gamma_1^2  \gamma_0^{-1} \Phi_m(\sigma_{\mathcal A}).
  \end{equation}
\end{corollary}
\begin{proof}
 The first inequality is clear since $\psi_{m}(\cdot,\mathcal A)$ belongs to the   family $\mathcal{G}_{m} (\Omega, \mathcal A)$ whose upper envelope is  $G_m (\cdot,\mathcal A,\Omega)$.

Let us prove the second one. Indeed, let $v\in\mathcal{G}_{m} (\Omega, \mathcal A)$ and let  $(a,\nu) \in \mathcal A$ be fixed. By definition there exists $C_v > 0$ such that for any $z \in \Omega$ 
 
 $$
 v(z)\leqslant \phi_m(z,\mathcal A) + C_v  \leq \nu \Phi_m(z-a) +C_v.
 $$ 
Fix  $0<\varepsilon < 1$ and choose $r_0 >0$ so that $\nu \Phi_m(r_0) + (1+\varepsilon) C_v = 0$. 
 Then  for any $0 < r < r_0$, we obtain
 $$
 (1+\varepsilon)v(z)\leqslant \nu \Phi_{m}(z-a) \, \, \, \text{on} \, \, \partial B(a,r).
 $$
   On the other hand,  fix   $ 0 < \delta < \theta_m(d(a),\nu^{-1})$. Then   $B (a,\theta_m(\delta ,\nu)) \Subset \Omega$ and $\nu \Phi_{m}(\cdot - a)  - \Phi_m(\delta) \geq 0$ on $ \bar \Omega \setminus B(a,\theta_m(\delta ,\nu))$, hence on $\partial \Omega$. 
   
  Therefore for any $z \in  \partial(\Omega\smallsetminus \B(a,r))$, we have 
 $$
 (1+\varepsilon)v(z)\leqslant  \nu  \, \Phi_{m}(z-a)- \Phi_m(\delta).
  $$
   Since $(1+\varepsilon)v,$ and $\nu \, \Phi_{m}(\cdot - a) - \Phi_m(\delta)$ are bounded $m$-subharmonic functions on $\Omega\smallsetminus \B(a,r)$ and $\nu \, \Phi_{m}(\cdot - a) - \Phi_m(\delta)$ is maximal on $\Omega \setminus B(a,r)$,  by applying  the comparison principle Proposition \ref{prop:CP1}, we get,
   $$
   (1+\varepsilon)v(z)\leqslant \nu \, \Phi_{m}(z-a) - \Phi_m(\delta) \, \, \,  \text{on} \, \, \,  \Omega \smallsetminus \bar{B}(a,r).
   $$
 Since $r>0$ is arbitrary small, it follows that 
 $$ 
 (1+\varepsilon)v(z)\leqslant \nu  \, \Phi_{m}(z-a) - \Phi_m(\delta),
 $$
  in $\Omega $. Since $\varepsilon >0$ is arbitrary, we conclude that $v(z)\leqslant \nu  \, \Phi_{m}(z-a) - \Phi_m(\delta)$ in $\Omega$.  Hence $G_m(z,\mathcal A,\Omega) \leq \phi_m(z,\mathcal A) - \Phi_m(\delta)$ in $\Omega$. 
 
Therefore for any $\delta \leq \delta_{\mathcal A}$, we have
$$
G_m (z,\mathcal A,\Omega)\leqslant  \,  \phi_{m}(z,\mathcal A) - \Phi_m(\delta), \, \, \, \text{in} \, \, \Omega.
$$
This implies the second inequality in (\ref{eq:FIneq4}). The inequality (\ref{eq:FIneq5})  follows from  (\ref{eq:FIneq4}) and (\ref{eq:FIneq3}). This proves the statement of the corollary.
  \end{proof}
     
  \smallskip
 \smallskip

   \subsection{Boundary behaviour of the Green function} 

Let us first recall a definition. A bounded open domain $\Omega \Subset \C^n$  is said to be $m$-hyperconvex ($1 \leq m \leq n$) if it admits a negative $m$-subarmonic exhaustion $\rho : \Omega \longrightarrow ]-\infty , 0[$. 
 
 Recall that for $\delta_0 > 0$ small enough and $0 < \gamma_0 <  \gamma_1$ fixed,
  $$
  \mathcal E (\delta_0,\gamma_0, \gamma_1) := \{ \mathcal A \subset \Omega \times \R^+ \, ; \, \delta (\mathcal A;\partial\Omega) \geq \delta_0, \inf_ {a \in A} \nu(a) \geq \gamma_0, \sum_{a \in A} \nu(a) \leq \gamma_1\}.
  $$  
  Observe that  $\delta (\mathcal A;\partial\Omega) \geq \delta_0$  iff $ A_{\delta_0} := \bigcup_{(a,\nu) \in \mathcal A} {B} (a,\theta_m(\delta_0, \nu) )\Subset \Omega$.
 Moreover if $\mathcal A \in \mathcal E (\delta_0,\gamma_0, \gamma_1)$, then $A$ is a finite with cardinality $p := \vert A\vert  \leq \gamma_1 \gamma_0^{-1}$.
 \begin{proposition}  \label{prop:BdValues} Assume that $\Omega \Subset \C^n$ is a bounded $m$-hyperconvex domain. Then we have 
 \begin{equation} \label{eq:BdValues}
 \lim_{z\to \partial \Omega} \left(\inf_{\mathcal A \in \mathcal E (\delta_0,\gamma_0, \gamma_1)} \, G_m(z,\mathcal A,\Omega)\right)  = 0.
 \end{equation}
\end{proposition}
\begin{proof}
 Recall that for $\delta > 0$ we have 
 $$
 A_\delta :=  \bigcup_{(a,\nu) \in \mathcal A} \B (a,\theta_m(\delta, \nu)) = \{z \in \C^n \, ; \, \phi_m(z,\mathcal A) < \Phi_m(\delta)\},
 $$ 
 and so if  $\mathcal A \in  \mathcal E (\delta_0,\gamma_0, \gamma_1) $ and $0 < \delta  <  \delta_0$, then $A_\delta \Subset \Omega$. 

Fix $0 < \delta_1 < \delta_0$ and observe that  for any $z \in \partial A_{\delta_1} \Subset \Omega$, we have $\phi_m(z,\mathcal A) = \Phi_m(\delta_1)$, hence $\psi_m (z,\mathcal  A) \geq  p  \Phi_m(\delta_1) \geq \gamma_1 \gamma_0^{-1}\Phi_m(\delta_1) $. 

 Let $\rho$ be a negative $m$-subharmonic defining function for $\Omega$. One can choose a large constant $C = C (\gamma_0,\gamma_1,\delta_1)  > 1$ so that  
 $  C \rho (z) \leq  \gamma_1 \gamma_0^{-1}\Phi_m(\delta_1)$ in $\partial A_{\delta_1} \Subset \Omega$ and then $C \rho \leq \psi_m (\cdot,\mathcal A)$ in $\partial A_{\delta_1}$.
 Then by the gluing principle, the  function defined by  
 \begin{equation} \label{eq:subGreen}
 v(z)=\left\{\begin{array}{lcl} {\psi}_m(z, \mathcal A) \ \   \text{ on}   \  A_{\delta_1} \\\\
  \max\big\{C\rho(z), {\psi}_m(z,\mathcal A)\big\}\ \  \text{ on}  \ \Omega\smallsetminus  A_{\delta_1}
  \end{array}\right.
\end{equation}
is a negative $m$-subharmonic function in  $\Omega$ which belongs to $\mathcal G_m (\Omega,\mathcal A)$  and then $v \leq G_m(\cdot,\mathcal A,\Omega) $ in $\Omega$.

Therefore for any $\mathcal A \in \mathcal E (\delta_0,\gamma_0, \gamma_1) $ we have  $C \rho \leq G_m(\cdot,\mathcal A,\Omega) \leq 0$ on $\Omega \setminus  A_{\delta_1}$, which proves the required property since $\lim_{z \to \partial \Omega} \rho (z) = 0$. 
 \end{proof}
 \smallskip
 \smallskip

\subsection{The Hessian measure of the Green function} 
Here we prove the following property.
\begin{proposition} \label{prop:Hessian}
Let $\Omega \Subset \C^n$ be a bounded domain and $\mathcal A \subset \Omega \times ]0,+\infty[$ a finite set of weighted poles. Then 
\begin{equation} \label{eq:Hessian}
(dd^c G_m(\cdot,\mathcal A,\Omega))^m \wedge \beta^{n-m} = c_{n,m}  \sum_{(a,\nu) \in \mathcal A} \nu^m \delta_{a}.
\end{equation}
\end{proposition}
\begin{proof}
We first show that the function $ G := G_m (\cdot,\mathcal A,\Omega)$ is a maximal $m$-subharmonic function on $\Omega\smallsetminus  A$. We proceed by the usual balayage process. Fix an euclidean ball $B \Subset \Omega\smallsetminus A$. 
Since $G $ is a bounded $m$-subharmonic function in a neighborhood of $\bar B$, we claim that there exist ${\widehat{G}} \in\mathcal{SH}_{m}(\Omega)\cap L^{\infty} (\Omega)$ such that $(dd^{c} \widehat{G})^{m}\wedge\beta^{n-m}=0$  in the sense of currents on $B$, $\widehat{G} \geq G$ in $\Omega$ and $\widehat{G} =G$ in $\Omega \setminus B$. 
This is a classical balayage trick  which goes back to Bedford and Taylor in the case $m=n$ (see \cite[Proposition 9.1]{BT82}).  
Indeed,  if $G$ is continuous in a neighborhood of $\bar B$,  we can use \cite[Theorem 3.7]{Bl05} to obtain  ${v} \in \mathcal{SH}_m (B) \cap L^{\infty} (B)$ such that 
$(dd^{c} {v})^{m}\wedge\beta^{n-m}=0$  in the sense of currents on $B$ and ${v} = G$ in $\partial B$. By the comparison principle we have ${v} \geq G$ in $B$. Hence the function defined by $\widehat{G} = v$ in $B$ and $\widehat{G}  = G$ in $\Omega \setminus B$ satisfies the requirements of the claim.   
In the general case, we approximate $G$ by a decreasing sequence of continuous $m$-subharmonic functions $ ({H}_j)_{j \in \N}$ in  $\Omega$ (see \cite{Lu15}). By the previous construction we obtain a  sequence $\widehat{H}_j$ of bounded $m$-subharmonic functions in $\Omega$ such that for any $j \in \N$,
$(dd^{c} \widehat{H}_j)^{m}\wedge\beta^{n-m}=0$  in the sense of currents on $B$, $\widehat{H}_j \geq H_j$ in $\Omega$ and $\widehat{H}_j  = H_j$ in $\Omega \setminus B$. By the comparison principle, the sequence $( \widehat{H}_j)$ is a decreasing sequence of functions in $ \mathcal{SH}_m (\Omega) \cap L^{\infty} (\Omega) $ which converges a.e. in $\Omega$ to a function $ \widehat{G} \in \mathcal{SH}_m (\Omega) \cap L^{\infty} (\Omega)$. It is clear that this function satisfies the required properties as we claimed (see \cite[Proposition 9.1]{BT82} for more details).
 
Now define the function  $u$ by  $u=\widehat{G}$ in $B$ and $u=G$ in $\Omega\smallsetminus B$. Then  $v\in\mathcal{G}(\Omega,\mathcal A)$. Hence $v\leqslant G $ which implies that $ \widehat{G} \leqslant G $ in $B$. This proves that $ \widehat{G} = G $ in $B$ and then  $(dd^c G)^{m}\wedge\beta^{n-m}=(dd^{c}\widehat{G})^{m}\wedge\beta^{n-m}=0$ in the sense of currents on $B$. 

  We next prove the formula (\ref{eq:Hessian}). Indeed, by the previous analysis 
  %by a standard result from measure theory we already know that
 the measure $(dd^c G)^m \wedge\beta^{n-m} $ have a finite support contained in the finite set $A$. Hence it is a finite combination of Dirac masses at the points in ${A}$. It is then enough to compute its  mass at each point $a \in A$. 
Fix  $a \in  A$ and observe   thanks to Corollary \ref{Cor:Fest3} that   $\lim_{z \to a} \frac{G (z)}{\Phi_m (z-a)} = \nu (a)$. 
We can then apply  Lemma \ref{lem:residualmass} to obtain the following formula
$$
\int_{\{a\}}(dd^c G)^m \wedge\beta^{n-m} = \nu(a)^m \int_{\{a\}}(dd^c \Phi_m (\cdot-a))^m\wedge\beta^{n-m} = c_{n,m} \nu(a)^m.
$$
This implies the formula (\ref{eq:Hessian}).
\end{proof}
\smallskip
\smallskip

\subsection{A generalized comparison principle } To prove the uniqueness in theorem A, we will need to prove a more general comparison principle which deals with singular  $m$-subharmonic functions in the class $\mathcal{SH}_m^b(\Omega)$.
 
  Since $\Phi_m$ is a radial fundamental solution of the Hessian equation  (\ref{eq:Fsol}), it follows from the comparison principle that  the singularity of any given  $m$-subharmonic function $u$  in $\Omega$ at any given point $a \in \Omega$ is at worst like $\nu \Phi_m(z-a)$ for some  constant $\nu \geq 0$. Indeed, applying  the comparison principle and taking into account the formula (\ref{eq:Fsol}), it is easy to see that the function $r \longmapsto \max_{\bar \B (a,r)} u $ is an increasing convex function of the variable $t := \Phi_m(r)$ for $0 < r <  d (a,\partial \Omega)$. Then the following limit exists :
   \begin{equation} \label{eq:Lelongnumber}
   \nu_m (u,a) : =\lim\limits_{r\rightarrow 0^+} \frac{\max_{\bar \B (a,r)} u}{\Phi_m (r)} \in [0, + \infty[, 
    \end{equation}
where $\bar{\B} (a,r)=\{z\in \Omega \, ; \, |z-a| \leq r\}$ is the euclidean ball. By convexity, for any $ 0 < r < r_0 < d (a,\partial \Omega)$, we have
\begin{equation} \label{eq:ConvexityIneq}
\max_{\bar \B (a,r)} u - \max_{\bar \B (a,r_0)} u \leq  \nu_m (u,a) (\Phi_m (z-a) - \Phi_m (r_0)).
\end{equation}
   This means that the real number $\nu_m(u,a)$  measures the weight of the singularity of $u$ at the point $a$.
   
  Lelong numbers  of  the $m$-positive current $dd^c u$ associated to a $m$-subharmonic function $u$  was introduced in \cite{WW16} and its relationship to the mean values of $u$ on spheres and balls was given in \cite{BG18}.

\smallskip

 We first prove the following elementary lemma.
 
  \begin{lemma} Let   $\mathcal A  := \{(a,\nu(a)) \, ; \, a \in A\} \subset \Omega \times \R^+$ be a finite weighted set.  Then 
  for any $a \in A$, $\nu_m (G_m(\cdot,\mathcal A,\Omega), a) = \nu (a)$.
  
   Moreover if $u \in \mathcal {SH}_m(\Omega)$,  $u \leq M$ in $\Omega$ and $\nu_m (u,a) \geq \nu (a)$ for any 
   $a \in A$,  we have $u \leq M + G_m(\cdot,\mathcal A, \Omega)$ in $\Omega$. In particular we have
   \begin{equation} \label{eq:LelongFormula}
   G_m(\cdot,\mathcal A, \Omega) = \sup \{ u \, ; \, u \in \mathcal{SH}_m^- (\Omega) ,\forall a \in A,  \nu_m (u,a) \geq \nu(a) \}.
   \end{equation}
   \end{lemma}
    \begin{proof} 
    
    Fix a point $(a,\nu(a)) \in \mathcal A$. Then by (\ref{eq:FIneq4}), we have for any $z \in \Omega$
    $$
    \psi_m(z,\mathcal A) = \nu(a) \Phi_m(z - a) + g (z) \leq G_m (z,\mathcal A,\Omega) \leq \nu(a) \Phi_m(z - a)  - \Phi_m(\delta_A),
    $$
    where $ g (z) := \sum_{b \in A, b \neq a}  \nu(b) \Phi_m(z - b) $ is a bounded $m$-sh function in a neighborhood of  the point $a$. 
    This implies that $\nu_m (G_m(\cdot,\mathcal A,\Omega), a) = \nu(a)$.
 
   Fix a point $(a,\nu(a)) \in \mathcal A$. By the convexity inequality (\ref{eq:ConvexityIneq}), it follows that  there exists constant $C_a, r_a > 0$ such that
 $$
 u (z) \leq \nu_m (u,a) \Phi_m(z-a) + C_a \leq \nu (a) \Phi_m(z-a) + C_a,
 $$
   for $z \in \bar \B (a,r_a) \subset \Omega$.  Since  $u \leq M$ in $\Omega$ and $ z \longmapsto \Phi_m(z-a) $ is bounded from below on $\Omega \setminus B(a,r_a)$, it follows that  there exists a constant $C'_a >0$ such that $ u(z)  \leq \nu(a) \Phi_m(z-a) + C'_a$ for any $z \in \bar \Omega$. Set $C' := \max_{a \in A} C'_a$. Then we have  $ u (z) \leq \phi_m (z,\mathcal A) + C'$ for any $ z \in \Omega$. Hence $u - M \in \mathcal G_m(\Omega,\mathcal A)$ 
  and then $u - M \leq G_m(\cdot,\mathcal A,\Omega)$ in $\Omega$. 
  
  The formula (\ref{eq:LelongFormula}) follows immediately from the above analysis.  \end{proof}
  
 We can easily prove the following Lemma (see \cite{Ze97}).
 \begin{lemma} Let  $u \in \mathcal {SH}^b_m(\Omega)$.  Then we have
 \begin{equation}  \label{eq:massineq}
 c_{m,n}  \sum_{a \in A_u} \nu_m (u,a)^m \leq \int_{S_u} (dd^c u)^m \wedge \beta^{n-m},
 \end{equation}
 where $A_u := \{ a \in \Omega ;  \nu_m (u,a) > 0\}$ and $S_u := \{a \in \Omega \, ; \, u(a) = - \infty\}$.
 \end{lemma}
  \begin{proof}  By definition, there exists a compact set $K \subset \Omega$ such that $u $ is bounded in $\Omega \setminus K$ and then $A_u \subset S_u \subset K$. Modifying $u$ near the boundary, we can assume that $u = 0$ in $\partial \Omega$. 
  Let $A \subset A_u$ be a finite set and  $\mathcal A := \{(a,\nu_m(u,a)) \, ; \, a \in A\}$. By the previous lemma, we have 
 $u \leq G_m (\cdot,\mathcal A)$ in $\Omega$. By the comparison principle Proposition \ref{prop:CP2}, we deduce that
 $$
 \int_A (dd^c G_m (\cdot,\mathcal A,\Omega)^m \wedge \beta^{n-m} \leq \int_{S_u} (dd^c u)^m \wedge \beta^{n-m},
 $$ 
 and the inequality (\ref{eq:massineq})  follows from (\ref{eq:MongeAmperemass}).
  \end{proof}
  
  Following (\cite{Ze97}), we can prove the comparison principle which will imply the uniqueness of the Green function stated in Theorem A.
 \begin{proposition}  \label{prop:Uniqueness} Let $E \subset \Omega$ be a compact subset of Lebesgue measure $0$.  Let $u, v  \in \mathcal{SH}_m(\Omega)  \cap L^{\infty}_{loc} (\Omega \setminus E)$ such that $\liminf_{z \rightarrow \zeta }(u(z)- v(z))\geq 0$. Assume that the following properties hold
 
 $(i) $ $\int_E (dd^c u)^m \wedge \beta^{n-m} =  c_{m,n}  \sum_{a \in A_u} \nu_m (u,a)^m$,
 
 $(ii)$ $(dd^c v)^m \wedge \beta^{n-m} \geq (dd^c u)^m \wedge \beta^{n-m}$, in the sense of measures on $\Omega \setminus E$,
 
 $(iii)$ $\nu_m (v,a) \geq \nu_m (u,a)$ for any $a \in E$.
 
 Then $u \geq v$ in $\Omega$.  
 \end{proposition}
 
 \subsection{Proof of Theorem A} Fix a finite weighted set $\mathcal A \subset \Omega \times ]0, + \infty[$ and set $G := G_m (\cdot, \mathcal A, \Omega)$. By  Corollary \ref{Cor:Fest3}, the upper semi-continuous regularization $G^*$ of $G$ satisfies the inequality
 $$
\psi_m (z, \mathcal A) \leq G^*(z) \leq \phi_m (z, \mathcal A) - \Phi_m (\delta_A),
 $$
 since $\phi_m$ is upper semi-continuous in $\Omega$. Hence $G^* \in \mathcal G_m (\Omega, \mathcal A)$ and then $G^* \leq G$ in $\Omega$. This implies that $G= G^*$ is $m$-subharmonic in $\Omega$ and satisfies the inequality (\ref{eq:global-estimate}).
  Proposition \ref{prop:BdValues} implies (\ref{eq:boundary-estimate}), the formula (\ref{eq:MongeAmperemass}) follows from Proposition \ref{prop:Hessian} and uniqueness follows from Proposition \ref{prop:Uniqueness}. This proves Theorem A.

 \section{Modulus of continuity of the Green function}
 
 The goal of this section is to prove Theorem B. This will be done in several steps. Let's outline the main steps.
 There are three steps starting from the following obvious inequality : for any $z, z' \in \Omega $ and $\mathcal A \subset \Omega \times \R^+$, $\mathcal A'  \subset \Omega \times \R^+$, we have 
 
 \begin{eqnarray} \label{eq:twoterms}
 \vert \exp G(z', \mathcal A') - \exp G(z, \mathcal A)\vert & \leq  & \vert \exp G(z', \mathcal A') - \exp G(z, \mathcal A')\vert  \nonumber \\
 &&  + \, \, \vert \exp G(z, \mathcal A') - \exp G(z, \mathcal A)\vert,
 \end{eqnarray}
 where $G(z,\mathcal A) := G_m (z,\mathcal A,\Omega)$.
 
 The {\it first step} done  in Section 4.1  consists in  estimating the  modulus of continuity of  the weight functions $\phi_m (z,\mathcal A)$  and   $\psi_m (z,\mathcal A)$ in terms of $z$ and $\mathcal A$ (see in Lemma \ref{lem:Flem}).
 
 The {\it second step} is done in Section 4.2. We use the first step to obtain an estimate of the first term on the right hand side of the inequality (\ref{eq:twoterms}) by a function  of $r := \vert z-z'\vert$ for $z,z' \in \Omega$, uniformly in $\mathcal A' \in \mathcal F (\delta_0, \gamma_0,\gamma_1, \sigma_0)$ (see Theorem \ref{thm:Holder-espace}). 
 
 This is the difficult step in the proof of Theorem B.  Here we use the classical technique of perturbation of the domain due to J.L. Walsh. This argument became classical and has been used originally in \cite{Wal68}  as well as in many other works to prove continuity of various envelopes (e.g.  \cite{Lel89}). 
 
 However since we want to get a precise control on the modulus of continuity of the Green function, we need to use an extra argument based on the subextension trick using Proposition \ref{prop:Subext}. 
 
The {\it third step} done in Section 4.3 is easier. We use the estimates proved in the first step  to estimate the second term on the right hand side of the inequality (\ref{eq:twoterms}) by a function of $r := d_H (\mathcal A,\mathcal A')$ for $\mathcal A' \in \mathcal F (\delta_0, \gamma_0,\gamma_1, \sigma_0)$, uniformly in $z \in \bar \Omega$ (see Theorem \ref{thm:Holder-espacepondere}). 
   \smallskip
   
As far as we know the idea of using  the maximal subextention argument in getting precise modulus of continuity of an envelope   is new and we believe it may be used in other contexts.  

 \subsection{Equicontinuity of the weighted functions}
 
 The first step in the proof of Theorem B will consist in proving Lemma \ref{lem:Flem} below.
 
 Fix $\mathcal A \subset \Omega \times \R^+$ and $\mathcal A'  \subset \Omega \times \R^+$ and recall that  
$$d_H ( \mathcal A , \mathcal A') \leq r 
\Longleftrightarrow \mathcal A \subset \mathcal V_r(\mathcal A') \, \text{and} \, \, \mathcal A' \subset \mathcal V_r(\mathcal A),
$$
where $\mathcal V_r(\mathcal A) := \bigcup_{x \in \mathcal A} \bar{\mathcal B} (x,r)$ is the $r$-neighborhood of $\mathcal A$ in $\Omega \times \R^+ \subset \C^n \times \R$.  Here $\bar{\mathcal B} (x,r)$ is the ball in $\Omega \times \R^+$ of center $x = (a,\nu) \in \Omega \times \R^+$ and radius $r$  for the distance $d_1 (x,x') :=  \vert a - a'\vert  + \vert \nu - \nu'\vert $, where $x = (a,\nu) \in \Omega \times \R^+$ and $x' := (a',\nu') \in \Omega \times \R^+$.

 In particular,  $\mathcal A \subset  \mathcal V_r(\mathcal A')$ iff  for any $(a,\nu) \in \mathcal  A$, there exists $(a',\nu') \in \mathcal A'$ such that $\vert a - a'\vert + \vert \nu - \nu'\vert \leq r$.    

   \smallskip
   For fixed $0 < \gamma_0 < \gamma_1$,  $\sigma_0 > 0$ and  $\delta_0 > 0$ small enough, we define two families of weighted sets.
   
   Recall that   $ \mathcal{E} (\delta_0,\gamma_0,\gamma_1)$ be the family of weighted sets  $ \mathcal A \subset \Omega \times \R^+,$ satisfying the following conditions:
   $$
   \delta (\mathcal A,\Omega) \geq \delta_0,   \, \, \inf_{(a,\nu) \in \mathcal A} \nu \geq \gamma_0,  \, \, \sum_{(a,\nu) \in \mathcal A} \nu \leq \gamma_1,
   $$
   and
   $$
   \mathcal{F} (\delta_0,\gamma_0,\gamma_1,\sigma_0) := \{\mathcal A  \, ;  \,  \mathcal A \in  \mathcal{E} (\delta_0,\gamma_0,\gamma_1), \sigma_{\mathcal A} \geq \sigma_0\}\cdot
   $$
    Recall the following definition 
 $$
 A_{\delta} := \bigcup_{(a,\nu) \in  \mathcal A}  B (a, \theta_m (\delta, \nu)) = \{z \in \C^n ; \phi_m(z,\mathcal A) < \Phi_m(\delta)\}\cdot
 $$
 We also define the following function
 \begin{equation}\label{eq:weigth-distance}
f_m (t) :=    \left\{\begin{array}{lcl} 
   t^{-2s - 1},  &\hbox{if}\  1 \leq m < n \, (s> 0),  \\
t^{- 1 \slash \gamma_0} & \hbox{if}\  m = n,
\end{array}\right.
\end{equation}

 We first prove the following lemma.
   \begin{lemma}  \label{lem:Flem} 
  Fix   $\mathcal A, \mathcal A' \in \mathcal E(\delta_0,\gamma_0,\gamma_1)$   and $0 < \delta \leq \delta_0$. Then  for any $z \notin A_{\delta}$, any $z' \notin A'_{\delta}$, we have 
 \begin{equation} \label{eq:FIneq6}
  \phi_m(z',\mathcal A') \leq \phi_m (z,\mathcal A) + L_m  \, f_m (\delta)   \left(\vert z - z'\vert + d_H( \mathcal A ;\mathcal  A')\right),
 \end{equation}
 and
 \begin{equation} \label{eq:FIneq7}
 \psi_m(z',\mathcal A') \leq \psi_m (z,\mathcal A)  +  L'_m    \, f_m (\delta) \left(\vert z - z'\vert + d_H( \mathcal A ;\mathcal  A')\right), 
 \end{equation}
where  $L_m  $ and $L'_m $ are uniform constants.
\end{lemma}
The constants $L_m$ and $L'_m$ are given respectively by (\ref{eq:Lm}) and (\ref{eq:L'm}) when $m < n$ and by (\ref{eq:Ln}) and (\ref{eq:L'n}) when $m=n$.
\begin{proof} Indeed,  fix $0 < \delta \leq \delta_0$ and  $z \notin A_{\delta}$. Then there exists $(a,\nu) \in \mathcal A$ such that
 $$
   \phi_m (z, \mathcal A) = \nu \Phi_m( z - a).
 $$
 By definition  there exits $(a',\nu') \in \mathcal A'$ such that 
 $\vert a - a'\vert + \vert \nu - \nu'\vert \leq  d_H (\mathcal A;\mathcal A')$.  Hence
$$
 \phi_m(z',\mathcal A') -  \phi_m (z, \mathcal A) \leq   \nu' \Phi_m(z' - a') -  \nu \Phi_m(z - a).
 $$
 Now observe that since $z' \notin A'_{\delta}$ and $z \notin A_{\delta} $, we have $\vert z- a \vert \geq  \theta_m(\delta, \nu) \geq \theta_m (\delta, \gamma_0)$ and $\vert z' - a'\vert \geq   \theta_m (\delta, \gamma_0) $. On the other hand, write 
 
  \begin{eqnarray*} 
  \nu' \Phi_m(z' - a') -  \nu \Phi_m(z - a)
 & = & \nu' (\Phi_m(z' - a') -  \Phi_m(z - a)) \\
 &+ & (\nu' - \nu) \Phi_m(z - a),
\end{eqnarray*}
and observe that  $\vert \Phi_m( t') - \Phi_m( t) \vert \leq   \Phi'_m(t_0)  \vert t'-t\vert$ for any real numbers $t, t' \geq t_0 :=  \theta_m (\delta, \gamma_0) > 0$. 

Then it follows that for  $z' \notin A'_{\delta}$ and $z \notin A_{\delta} $, we have
\begin{eqnarray} \label{eq:Ineqcomparison}
  \nu' \Phi_m(z' - a') -  \nu \Phi_m(z - a) & \leq  & \nu' \Phi_m' \circ \theta_m (\delta,\gamma_0) \left( \vert z' - z\vert  + \vert a'-a\vert\right) \nonumber \\
  && - \Phi_m \circ \theta_m (\delta,\gamma_0)   \vert \nu' - \nu\vert.
 \end{eqnarray}

We proceed to prove the estimate (\ref{eq:FIneq6}), by considering the two cases separately.

1.  {\it The case  $1 \leq m < n$. } In this case $\theta_m (\delta,\gamma_0) = \gamma_0^{1\slash 2s} \delta$ and $\Phi_m'(t) =  2s t^{-2s-1}$. Hence the equation (\ref{eq:Ineqcomparison}) yields for  $z' \notin A'_{\delta}$ and $z \notin A_{\delta} $,
\begin{eqnarray*} 
  \nu' \Phi_m(z' - a') -  \nu \Phi_m(z - a) & \leq  & 2s \gamma_1 \gamma_0^{- (2s + 1)\slash 2 s} \delta^{-2s -1} \left( \vert z' - z\vert  + \vert a'-a\vert\right) \nonumber \\
  && +  \gamma_0^{-1} \delta^{-2s} \vert \nu' - \nu\vert.
 \end{eqnarray*}
This implies that  for any $0 < \delta \leq \delta_0,$ any $z \notin A_{\delta} $ and $z' \notin A'_{\delta}$,  we have
 $$
 \phi_m(z',\mathcal A') -  \phi_m (z, \mathcal A) \leq L_m  \, \delta^{- 2s - 1}  \left(\vert z' - z\vert +  d_H (\mathcal A;\mathcal A')\right),
 $$
 where 
 \begin{equation} \label{eq:Lm}
 L_m := \max \{ 2s \gamma_1 \gamma_0^{-(2s+1)\slash 2s} ,\delta_0  \gamma_0^{-1}\}.
 \end{equation} 
 This proves the estimate (\ref{eq:FIneq6}) when $m < n$. 
 
 \smallskip
 
 2.  {\it The case  $ m =n$. }  In this case $\theta_n (\delta,\gamma_0) = R_0(\delta\slash R_0)^{1 \slash \gamma_0}$ and $\Phi'_n(t) = R_0 \slash  t$. Hence the equation (\ref{eq:Ineqcomparison}) yields  for  $z' \notin A'_{\delta}$ and $z \notin A_{\delta} $, we have
 \begin{eqnarray} \label{eq:Ineqcomparison2}
  \nu' \Phi_n(z' - a') -  \nu \Phi_n(z - a) & \leq  & \gamma_1 (R_0\slash \delta)^{1 \slash \gamma_0} \left( \vert z' - z\vert  + \vert a'-a\vert\right) \nonumber \\
  && + \, \gamma_0^{-1}  \log (R_0 \slash \delta)  \vert \nu' - \nu\vert.
 \end{eqnarray}
 Since $\gamma_0^{-1}  \log (R_0 \slash \delta) \leq  (R_0\slash \delta)^{1 \slash \gamma_0}$ for $0 < \delta < R_0$, it follows that 
 for any $0 < \delta \leq \delta_0,$ any $z \notin A_{\delta} $ and $z' \notin A'_{\delta}$,  we have
 $$
 \phi_n (z',\mathcal A') - \phi_n (z,\mathcal A) \leq  \gamma_1 (R_0\slash \delta)^{1 \slash \gamma_0} \left( \vert z' - z\vert  + d_H(\mathcal A,\mathcal A')\right),
 $$
 which proves the estimate (\ref{eq:FIneq6}) with the constant
 \begin{equation} \label{eq:Ln}
 L_n := \gamma_1  R_0^{1 \slash \gamma_0} 
 \end{equation}
 
 Now we prove the estimate (\ref{eq:FIneq7})  in the same way. Indeed, 
  observe that by definition  for any $(a,\nu(a)) \in \mathcal A$ there exists $(b(a),\mu(a)) \in \mathcal A'$ such that 
 $\vert a - b(a)\vert + \vert \nu- \mu(a)\vert \leq  d_H (\mathcal A;\mathcal A')$.  
 
 For $z \notin (A \cup A')$ we have
 $$
 \psi_m(z',\mathcal A') =  \sum_{(a',\nu')  \in \mathcal A'} \nu' \Phi_m(z' - a')  \leq \sum_{a \in A}  \mu (a) \Phi_m(z' - b(a)). 
 $$ 
 Then using (\ref{eq:Ineqcomparison}), we get
\begin{eqnarray*}
 \psi_m(z',\mathcal A') -  \psi_m (z, \mathcal A) & \leq&  \gamma_1 \Phi_m' \circ \theta_m(\delta,\gamma_0) \left( \vert z' - z\vert  + d_H(\mathcal A,\mathcal A')\right)  \\
  && - \gamma_1 \gamma_0^{-1}\Phi_m\circ \theta_m (\delta,\gamma_0) d_H(\mathcal A,\mathcal A').
\end{eqnarray*} 

1. Assume first that $1 \leq m < n$. Then   it follows as before that for any $0 < \delta \leq \delta_0$, any $z \notin A_{\delta} $ and $z' \notin A'_{\delta}$ we have
$$
 \psi_m(z',\mathcal A') -  \psi_m (z, \mathcal A) \leq L'_m  \, \delta^{-2s-1} \left(\vert z' - z\vert +  d_H (\mathcal A;\mathcal A')\right),
$$
where 
\begin{equation} \label{eq:L'm}
L'_m :=   \max \{ 2s \gamma_1 \gamma_0^{-(2s+1)\slash 2s},   \gamma_1  \gamma_0^{-2} \delta_0\}.
\end{equation}

2. Assume now that $m = n$.  The same computation shows that 
for any $0 < \delta \leq \delta_0$, any $z \notin A_{\delta} $ and $z' \notin A'_{\delta}$ we have
$$
 \psi_n(z',\mathcal A') -  \psi_n (z, \mathcal A) \leq L'_n \delta^{- 1 \slash \gamma_0} \left(\vert z' - z\vert +  
 d_H (\mathcal A;\mathcal A')\right),
$$
where 
\begin{equation} \label{eq:L'n}
L'_n :=   \gamma_1 R_0^{1 \slash \gamma_0}.
\end{equation}

This proves the estimate (\ref{eq:FIneq7}).
\end{proof}
\begin{remark} The previous result shows that for any open subset $D \Subset \Omega$, the weighted function $\phi_m$ and $\psi_m$ are Lipschitz continuous in $(\bar \Omega \setminus D) \times \mathcal F_D$, where $\mathcal F_D = \mathcal F_D(\delta_0,\gamma_0,\gamma_1,\sigma_0)$ is the family of  sets $\mathcal A \in \mathcal F (\delta_0,\gamma_0,\gamma_1,\sigma_0)$ such that $\mathcal A \subset D \times ]0,+\infty[$ endowed with the Hausdorff distance.
\end{remark}   

\smallskip

  \subsection{Equicontinuity  in the space variable}
  The second step in the proof of Theorem B will consist in proving the following result.
  
   \begin{theorem} \label{thm:Holder-espace}  Assume that $\Omega \Subset \C^n$ is a bounded  $m$-hyperconvex domain of Lipschitz type. Let $0 < \gamma_0 < \gamma_1$, $\delta_0 > 0$, and $\sigma_0 > 0$ be fixed. Then the family $\{\exp G_m (\cdot,\mathcal A, \Omega) \, ; \, \mathcal A \in \mathcal F (\delta_0,\gamma_0,\gamma_1,\sigma_0)\}$ is equicontinuous in $\bar \Omega$. More precisely we have :
   
   1. If $1 \leq m \leq n-1$,  for any $0 < \tau < \frac{2s}{2s+1}  $, there exists constants $L > 0$ and $r_1 > 0$, depending on $ (m,n,\tau,\delta_0,\gamma_0,\gamma_1)$ such that for any $\mathcal A \in \mathcal F(\delta_0,\gamma_0,\gamma_1,\sigma_0)$, and any $z, z' \in \bar \Omega$  with $\vert z - z'\vert \leq r \leq r_1$, we have
   \begin{equation} \label{eq:Holder-espace}
\vert  \exp G_m(z',\mathcal A, \Omega) - \exp G_m (z,\mathcal A,\Omega) \vert \, \leq \, L \cdot r^{\tau},
\end{equation}
 
2. If $m=n$, for any $0 < \alpha < 1$,   there exists a constant  $L > 0$ and $r_1 > 0$, depending on $(m,n,\alpha, \delta_0,\gamma_0,\gamma_1)$  such that for any $z, z' \in \Omega$ with $\vert z - z'\vert \leq r \leq r_1$, 
\begin{equation} \label{eq:equicontespace}
\vert  \exp G_n (z',\mathcal A, \Omega) - \exp G_n (z,\mathcal A,\Omega) \vert \, \leq \, L \cdot  (\log (R_1 \slash r))^{- \alpha},
\end{equation}
where $R_1 := R_0^{1 \slash \gamma_0}$.
\end{theorem}

\begin{proof}
 As we allready said in the beginning of section 4, we will use the classical technique of perturbation of domains due to J.B. Walsh \cite{Wal68}. 
 
Fix  $\zeta \in \C^n$ with $\vert \zeta\vert$ small enough, $\mathcal A \in \mathcal F (\delta_0,\gamma_0,\gamma_1,\sigma_0)$ and set for simplicity $G (z) := G_m (z,\mathcal A,\Omega)$ for $z \in \Omega$. 

 Define the perturbed domain $\Omega^\zeta := \{ z \in \Omega ; z + \zeta \in \Omega\}$ and set $G^\zeta (z) := G (z + \zeta)$ for $z \in \Omega^\zeta$. Then $G^\zeta$ is $m$-subharmonic in $\Omega^\zeta$. 
 
 The idea is to modify suitably $G^\zeta$  to produce a $m$-subharmonic function in the class $\mathcal G_m (\Omega,\mathcal A)$ which enables  to compare $G^\zeta$  and  $G$ in $\Omega \cap \Omega^\zeta$. 
 
 To this end we need to construct a $m$-subharmonic function in $\Omega$ close to the perturbed function $G^\zeta$  in the domain $\Omega \cap \Omega^\zeta$ and having the same singularities as $G$.   This will be done in two steps. First by  a max construction we produce such a function  in the domain $\Omega \cap \Omega^\zeta$. The main difficulty is to "extend" this function to $\Omega$. This is  based on a tricky argument using the maximal subextention method given in  Proposition \ref{prop:Subext} that we will explain below.

\smallskip

As in the previous proofs, we consider two cases.

\smallskip

  1. {\it The case $1 \leq m<n$.} We first assume that $ 1 \leq m < n$. By (\ref{eq:FIneq5}),  for any $z \in \Omega^\zeta$, we have
 $$
 G^\zeta (z) = G (z+ \zeta) \leq \psi_m(z + \zeta,\mathcal A) + \delta_0^{-2 s} + \gamma_1^2 \gamma_0^{-1} \sigma_0^{-2s}.
 $$
 To get rid of the constant in the right hand side, we introduce a small parameter $\varepsilon >0$. Fix $\varepsilon_0 > 0$ to be chosen later and let $0 < \varepsilon \leq \varepsilon_0$, $0 < \delta \leq \delta_0$ and observe that if $0 < r \leq \gamma_0^{1 \slash 2s} \delta \slash 2 $ and $\vert \zeta \vert \leq r$, then $A \subset \Omega^\zeta$ and for any $z \in \Omega^\zeta \setminus A_{\delta}$,  $z+ \zeta \notin A_{\delta\slash 2}$. 
By (\ref{eq:FIneq3}) and Lemma \ref{lem:Flem},  this implies that for any $z \in \Omega^\zeta \setminus A_{\delta}$, we have
  \begin{eqnarray*}
(1+ \varepsilon) G (z + \zeta) &\leq & (1+ \varepsilon) \psi_m(z + \zeta, \mathcal A) + (1+ \varepsilon) (\delta_0^{-2 s} + \gamma^2_1 \gamma_0^{-1} \sigma_0^{-2s})\\
 &\leq & (1+ \varepsilon) \psi_m (z, \mathcal A)  +  (1+ \varepsilon)  2^{2 s +1} L'_m \, r  \delta^{-2 s-1}   \\
 & + & (1+ \varepsilon) (\delta_0^{-2 s} + \gamma_1^2 \gamma_0^{-1} \sigma_0^{-2s}).
 \end{eqnarray*}

 Hence for  $0 < \varepsilon < \varepsilon_0$ and $z \in \Omega^\zeta \setminus A_{\delta}$,
 \begin{eqnarray*}
(1+ \varepsilon) G (z + \zeta)  &\leq & \psi_m (z, \mathcal A) + \varepsilon \psi_m(z,\mathcal A)  + C_0 +C_1 r  \delta^{-2 s-1},
   \end{eqnarray*}
  where $C_0 :=  (1+\varepsilon_0)  (\delta_0^{-2 s}  + \gamma_1^2 \gamma_0^{-1} \sigma_0^{-2s})$ and $C_1 :=  (1+\varepsilon_0)  2^{ 2s +1} L'_m$.

  Recall that $A_\delta = \{\phi_m(z,\mathcal A) < -\delta^{-2 s}\}$ and since $  \psi_m (z, \mathcal A) \leq \phi_m(z,\mathcal A)$, it follows that  for $z \in \Omega^\zeta \cap \partial A_{\delta}$, we have
  $$
  (1+ \varepsilon) G (z + \zeta)  \leq   \psi_m (z, \mathcal A) -\varepsilon \delta^{-2s} + C_0   + C_1 r \delta^{-2s -1},
  $$
    
 Set  $ r_0 := \min \{2^{- (2 s +1)\slash 2 s} \gamma_0^{(2 s +1)\slash 4 s^2}, \delta_0^{2s + 1}\}$ and  define for $0 < r\leq r_0$,  $ \delta  := r^{1 \slash (2s +1)}$ i.e. $r := \delta^{2 s + 1}$. Then for $0 < r \leq r_0$ we have  $r \leq \gamma_0^{1 \slash 2 s} \delta \slash 2$ and $\delta \leq \delta_0$. Now we define $ \varepsilon = \varepsilon (r) > 0$ so that $-\varepsilon \delta^{-2s} + C_0   + C_1 r \delta^{-2s -1} = 0$ i.e. 
 $$
 \varepsilon (r) :=  C_0  \delta^{2 s} +  C_1  \delta^{2 s}  = C_2 r^{2 s \slash(2s +1)},
 $$
 where $C_2 := C_0 + C_1$, so that  $ - \varepsilon \delta^{-2 s} +  C_0 +C_1 = 0$.
Since for $0 < r \leq r_0$, we have $\varepsilon (r) \leq C_2  r_0^{2s  \slash (2 s+1)}$, we choose $\varepsilon_0 := C_2   r_0^{2s  \slash (2 s+1)}$. 

Fix $0 < r \leq r_0$ and $\vert \zeta \vert = r$. Then  for any  $z \in \Omega^\zeta \cap \partial A_{\delta}$,  we have
  $$
  (1+ \varepsilon) G^\zeta (z)   \leq   \psi_m (z, \mathcal A).
  $$
  It is easy to see that one can find $0 < r_1 < r_0$ such that  if  $\vert \zeta \vert =r \leq r_1$  then $\zeta + A_\delta \subset A_{\delta_0} = \bigcup_{(a,\nu) \in \mathcal A} \B(a,\nu^{1\slash2s} \delta_0) \Subset \Omega$.  This implies that 
$A_\delta \Subset \Omega^\zeta$.
 Therefore the following function
\begin{equation} \label{eq:vzeta}
 v^\zeta (z)  :=   \left\{\begin{array}{lcl} 
  \psi_m (z,\mathcal A)  ,  &\hbox{in}\   A_\delta \\
  \max \{(1 + \varepsilon) G^\zeta(z) , \psi_m (z, \mathcal A)\} & \hbox{in}\  \Omega^\zeta \setminus A_\delta
\end{array}\right.
\end{equation}
is a negative $m$-subharmonic in $\Omega \cap \Omega^\zeta$ which has the right singularities on $A \subset  A_\delta \subset \Omega \cap \Omega^\zeta $.  

We would like to show that $v^\zeta \leq G + O (\vert \zeta\vert)$ in $\Omega \cap \Omega^\zeta$. This will be the case if we could extend $v^\zeta$ as a negative $m$-subharmonic function in $\Omega$. This is not clear but instead we can consider its maximal $m$-sh subextension to $\Omega$ defined as follows:
 \begin{equation} \label{eq:wzeta}
 w^\zeta := \sup \{ u \in \mathcal{SH}^-_m (\Omega) \, ; \, u \leq v^\zeta, \, \, \text{in} \, \,  \Omega \cap \Omega^\zeta \}.
 \end{equation} 
Since  $\psi_m (\cdot,\mathcal A) \leq v^\zeta$ in $\Omega^\zeta$, it follows that $w^\zeta$ is a well defined negative $m$-sh function in $\Omega$ which satisfies the inequalities $ \psi_m (\cdot,\mathcal A) \leq w^\zeta \leq v^\zeta$ in $\Omega \cap \Omega^\zeta$. Moreover since $ v^\zeta = \psi_m (\cdot,\mathcal A) $ in  $A_\delta$, it follows that    $w^\zeta =\psi_m (\cdot,\mathcal A)$ in $A_\delta$. This implies that $w^\zeta \leq G$ in $\Omega$. 

The goal is now to compare $v^\zeta$ and $w^\zeta$ in $\Omega \cap \Omega^\zeta$. We  first compare them on the boundary of $\Omega \cap \Omega^\zeta$ by finding a suitable sub-extension. Indeed let $\rho$ be a negative  $m$-subharmonic exhaustion function for $\Omega$ which is  Lipschitz continuous in $\bar \Omega$.  

Observe that by (\ref{eq:FIneq1}), for  any $z \notin A_{\delta_0}$ we have 
\begin{eqnarray*}
\psi_m(z, \mathcal A) &\geq & \phi_m(z,\mathcal A) +  \gamma_1 \gamma_0^{-1} \Phi_m(\delta_0) \geq - 2 \gamma_1 \gamma_0^{-1} \delta_0^{-2s}
\end{eqnarray*}
and choose a constant  $C = C(\delta_0,\gamma_0,\gamma_1) > 0$ such that 
$$
C \rho (z) \leq - 2 \gamma_1 \gamma_0^{-1} \delta_0^{-2s} \leq \psi_m (z,\mathcal A), \, \, \text{in} \, \, \partial  A_{\delta_0}.
$$
Therefore  the following function
\begin{equation} \label{eq:subObstacle}
 \theta (z)=\left\{\begin{array}{lcl} {\psi}_m(z,A) \ \   \text{ on}   \  A_{\delta_0} \\\\
  \max\big\{C \rho(z), {\psi}_m(z,\mathcal A)\big\}\ \  \text{ on}  \ \Omega\smallsetminus  A_{\delta_0}
  \end{array}\right.
\end{equation}
is a negative $m$-subharmonic function in  $\Omega$ which satisfies  $\theta \leq G(\cdot,\mathcal A) $ in $\Omega$. Since $\rho$ is a Lipschitz continuous, there exists a constant $M >0$ such that for any  $z, z'  \in  \bar \Omega$, we have  $\vert \rho (z) - \rho(z')  \vert \leq  M \vert z'-z\vert$. 
Hence 
$$
 C \rho (z) -   C M r \leq \theta (z + \zeta) \leq G(z+\zeta) =  G^{\zeta} (z),
$$ 
for $\vert \zeta \vert \leq r$ and $z \in \Omega \cap \Omega^\zeta \setminus A_{\delta_0}$.

Since $A_\delta \subset A_{\delta_0}$, it follows that the function $u  := (1+ \varepsilon)  \theta  - (1+\varepsilon) C M r$ is $m$-subharmonic in $\Omega$ and satisfies $u  \leq v^\zeta$ in $\Omega \cap \Omega^\zeta $.  Hence $u \leq w^\zeta$ in $\Omega$.

Recall that  if  $\vert \zeta \vert =r \leq r_1 \leq r_0$ then $r \leq \gamma_0^{1 \slash 2s}  \delta \slash 2$ with $\delta \leq \delta_0$, and $A_{\delta} \subset \Omega \cap \Omega^\zeta$. 
Observe that $\theta \geq 0$ in $ \Omega^\zeta \cap \partial \Omega$ and for $z \in \Omega \cap \partial \Omega^\zeta$,
$$
(1+\varepsilon) \theta (z) \geq (1+\varepsilon) C \rho (z) = (1+\varepsilon) C (\rho (z) - \rho(z+\zeta)) \geq  - (1+\varepsilon) C M r,
$$
Hence $(1+\varepsilon) \theta \geq - (1+\varepsilon) C M r$ in $\partial (\Omega \cap \Omega^\zeta)$ and then 
$$
w^\zeta \geq u =  (1+\varepsilon) \theta - (1+\varepsilon) C M r \geq  - 2 (1+\varepsilon) CM r \, \, \text{in} \, \, \partial (\Omega \cap \Omega^\zeta).
$$ 
Since $v^\zeta \leq 0$ in $\Omega \cap \Omega^\zeta$, it follows that
$$
v^\zeta  -  2 (1+ \varepsilon) CM r\leq w^\zeta,\, \, \text{in} \, \, \partial (\Omega \cap \Omega^\zeta).
$$
 Since $v^\zeta = \psi_m(\cdot,\mathcal A) \leq \theta \leq w^\zeta $ in $A_\delta$, we conclude that $v^\zeta -  2 (1+\varepsilon) CM r\leq w^\zeta $ in $\partial D$, where $D := (\Omega \cap \Omega^\zeta) \setminus \bar{A_\delta}$.

Now observe that obviously
$$v^\zeta -  2 (1+ \varepsilon) CM r\leq v^\zeta = w^\zeta , \, \, \text{in} \, \, \, D \cap \mathcal Q,
$$
 where  $\mathcal Q := \{w^\zeta = v^\zeta\}.$  Since by Proposition \ref{prop:Subext}, the Hessian measure $\mu_\zeta := (dd^c w^\zeta)^m \wedge \beta^{n-m}$ is carried by the contact set $\mathcal Q$, it follows that $v^\zeta -  2 (1+ \varepsilon) C M r\leq w^\zeta $,  $\mu^\zeta$-almost everywhere in $D$. By the domination principle, it follows that   $v^\zeta -  2 (1+ \varepsilon) C M r\leq w^\zeta $ in $ D$. Since $v^\zeta = w^\zeta = \psi_m (\cdot,\mathcal A)$ in $A_\delta$, it follows that $v^\zeta -  2 (1+ \varepsilon) C M r\leq w^\zeta $ in $\Omega \cap \Omega^\zeta$.
 
 This implies that $v^\zeta -  2 (1+ \varepsilon) C M r\leq G$ in $\Omega \cap \Omega^\zeta$, hence we obtain the following basic inequality
 \begin{equation} \label{eq:basic1}
 (1+ \varepsilon) G^\zeta \leq G + 2 (1+ \varepsilon) C M r,
 \end{equation}
  in $\Omega \cap \Omega^\zeta \setminus A_\delta$. 
 
Therefore  for any   $\vert \zeta \vert \leq r \leq r_1$ and for any $z \in (\Omega^\zeta \cap \Omega) \setminus A_\delta$,   we get 
\begin{eqnarray} \label{eq:Equ1}
G(z+ \zeta)  - G(z) & \leq &  -  \varepsilon  G(z + \zeta) + 2 (1+\varepsilon) C M  r \nonumber \\
&= & -  C_2 r^{2 s \slash(2s +1)} G(z + \zeta) + 2 (1+\varepsilon) C M  r.
\end{eqnarray}
 
 On the other hand, fix $\eta > 0$ and estimate  $-  G(z + \zeta) $ from above when  $z \notin A_{2 \eta}$. Indeed for such $z$ we have $ z' := z + \zeta \notin A_{\eta}$
if  $\vert \zeta \vert = r \leq \gamma_0^{1 \slash 2s} \eta$. 
 Then for any $a  \in A$ we have $\vert z+ \zeta  - a\vert \geq \vert z  - a\vert - r \geq \gamma_0^{1 \slash 2 s} \eta$.
 
 It follows from (\ref{eq:FIneq4}) that
  $$
  -G (z+\zeta) \leq  - \sum_{a \in A} \nu(a) \Phi_m(z + \zeta - a) \leq  \gamma_0^{-1} \gamma_1 \eta^{-2 s}.
   $$
 Recall that $\delta = r^{1 \slash (2 s + 1)}$, fix $\alpha >0$ small enough and set $\eta := \delta^\alpha \slash 2$. Choose  
 $0< r_2 < r_1$ such that for $r\leq r_2$ we have $r \leq \gamma_0^{1 \slash 2s} \delta^\alpha\slash 2$. Then  by the previous inequality for $\vert \zeta\vert \leq r \leq r_2$ and $z \in \Omega \setminus A_{\delta^\alpha}$, we have
 \begin{equation} \label{eq:Equ2}
  -G (z+\zeta) \leq \gamma_0^{-1} \gamma_1 \eta^{-2 s}.
 \end{equation}
 Since $A_ \delta \subset A_{\delta^\alpha}$, we have  $\Omega \setminus A_{\delta^\alpha} \subset \Omega \setminus A_{\delta}$. Then we can combine the estimates (\ref{eq:Equ1}) and (\ref{eq:Equ2}) to obtain the following estimate:
 for any $z \in \Omega \setminus A_{\delta^\alpha}$ and  $\vert \zeta \vert = r \leq r_2$ and $z \in \Omega^\zeta \setminus A_{\delta^\alpha}$, 
 \begin{eqnarray}
G(z+ \zeta) - G(z) &\leq &  C_2 r^{2 s \slash(2s +1)}   2^{2 s}   \gamma_1 \gamma_0^{-1} \, r^{-2 s \alpha \slash(2s+1)} \nonumber \\
&+&   2 (1+\varepsilon) C M r  \nonumber \\
& \leq &   C_3 r^{\tau(\alpha)},
\end{eqnarray}
where $\tau (\alpha) := (1-\alpha) \frac{2s}{2s+1}$ and 
$$C_3 :=  2^{2 s}   \gamma_1 \gamma_0^{-1} C_2 +   2 (1+ \varepsilon_0) C M r_1^{1-\tau(\alpha)}.
$$

Now we proceed to the proof of the uniform continuity of $\exp G_m (z,\mathcal A,\Omega)$ in $\bar \Omega \times \mathcal{F} (\delta_0,\gamma_0,\gamma_1,\sigma_0)$ with a control of its modulus of continuity.

Observe that for $x, y \in ]-\infty,0]$ we have $\vert e^x - e^y \vert \leq \vert x - y\vert$. Hence  for any $\mathcal{A}  \in \mathcal{F} (\delta_0,\gamma_0,\gamma_1,\sigma_0)$, $z \in \Omega \setminus A_{\delta^\alpha}$ and $z' \in \Omega \setminus A_{\delta^\alpha}$  such that  $\vert z'-z\vert \leq r  \leq r_2$, we have 
\begin{equation}\label{eq:Holder-exponent1}
 \exp G_m (z', \mathcal A,\Omega) - \exp G_m (z, \mathcal A,\Omega)  \leq L_1 r^{\tau},
\end{equation} 
where $r := \delta^{2 s +1}$ and $\tau := \tau (\alpha) := (1-\alpha) \frac{2s}{2s+1}$.

Now we want to estimate the left hand side of the inequality (\ref{eq:Holder-exponent1}) for $z \in A_{\delta^\alpha}$, $z' \in \Omega$  such that  $\vert z'-z\vert \leq r \leq r_2$.
Indeed, fix $ \mathcal A  \in \mathcal F$ and let $z \in A_{\delta^\alpha}$, $z' \in \Omega$  such that  $\vert z'-z\vert = r \leq r_2$.
Indeed there exists $(a,\nu) \in \mathcal A$ such that $\vert z - a\vert \leq \nu^{1\slash 2s} \delta^\alpha$ and then
 $$ 
 \vert z'-a\vert \leq \gamma_1^{1\slash 2s} r^{\alpha \slash (2s +1)} + r =: h(r).
 $$
Then by  (\ref{eq:FIneq4}), we have $G_m (z',\mathcal A,\Omega) \leq \frac{- \nu}{\vert z'-a\vert^{2s}} + \delta_0^{-2s}  \leq - \gamma_0 h(r)^{-2s} + \delta_0^{-2s}$ for any $z\in \Omega$. Hence
$$
 \exp G_m (z', \mathcal A, \Omega)  -  \exp G_m (z, \mathcal A, \Omega)  \leq \exp \left( -  \gamma_0 h(r)^{-2s} + \delta_0^{-2s}\right)
$$
Since $h(r) \simeq_{r \to 0^+}  r^{\alpha \slash (2s +1)}$, it follows that  there exists a constant $L_2 = L_2 (m,n,\alpha,\gamma_0,\delta_0)>0$ such that $\exp \left( -  \gamma_0 h(r)^{-2s} + \delta_0^{-2s}\right) \leq L_2 r^{\tau}$. Hence for any $z \in A_{\delta^\alpha}$, $z' \in \Omega$  such that  $\vert z'-z\vert = r \leq r_2$, we have
 \begin{equation}\label{eq:Holder-exponent2}
 \exp G_m (z', \mathcal A,\Omega)  -  \exp G_m (z, \mathcal A,\Omega)  \leq L_2 r^{\tau}.
\end{equation} 
 The inequality (\ref{eq:Holder-espace}) of the theorem follows from (\ref{eq:Holder-exponent1}) and (\ref{eq:Holder-exponent2}).
 
 \smallskip
 
 2. {\it The case $m=n$}. We  proceed in the same way. We again denote $G := G_n(\cdot,\mathcal A,\Omega)$. By (\ref{eq:FIneq5}),  for any $z \in \Omega^\zeta$, we have
 $$
 G^\zeta (z) = G (z+ \zeta) \leq \psi_n (z + \zeta,\mathcal A) + C_0,
 $$
 where $C_0 :=  \log (R_0\slash \delta_0) + \gamma_1^2 \gamma_0^{-1}  \log (R_0\slash \sigma_0)$. 
 
 To get rid of the constant in the right hand side, we introduce a small parameter $\varepsilon >0$. Fix $\varepsilon_0 > 0$ to be chosen later and let $0 < \varepsilon \leq \varepsilon_0$, $0 < \delta \leq \delta_0$ and observe that there exists $r_0 = r_0 (\gamma_0,R_0) > 0$ such that if $0 < r \leq r_0 \delta $ and $\vert \zeta \vert \leq r$, then $A \subset \Omega^\zeta$ and for any $z \in \Omega^\zeta \setminus A_{\delta}$,  we have $z+ \zeta \notin A_{\delta\slash 2}$. 
By (\ref{eq:FIneq3}) and Lemma \ref{lem:Flem},  this implies that for any $z \in \Omega^\zeta \setminus A_{\delta}$, we have
  \begin{eqnarray*}
(1+ \varepsilon) G (z + \zeta) &\leq & (1+ \varepsilon) \psi_n (z + \zeta, \mathcal A) + (1+ \varepsilon) C_0 \\
 &\leq & (1+ \varepsilon) (\psi_n (z, \mathcal A)  +  (1+ \varepsilon)  L'_n \, r  \delta^{-1\slash \gamma_0}   \\
 & + & (1+ \varepsilon) C_0.
 \end{eqnarray*}

 Hence for  $0 < \varepsilon < \varepsilon_0$ and $z \in \Omega^\zeta \setminus A_{\delta}$,
 \begin{eqnarray*}
(1+ \varepsilon) G (z + \zeta)  &\leq & \psi_n (z, \mathcal A) + \varepsilon \psi_n (z,\mathcal A)  + C'_0 +C'_1 r  \delta^{-1\slash \gamma_0},
   \end{eqnarray*}
  where $C'_0 :=  (1+\varepsilon_0) C_0$ and $C'_1 :=  (1+\varepsilon_0)  L'_n$.

  Recall that $A_\delta = \{\phi_m(z,\mathcal A) < \log (\delta\slash R_0)\}$ and since $  \psi_n (z, \mathcal A) \leq \phi_n (z,\mathcal A)$, it follows that  for $z \in \Omega^\zeta \cap \partial A_{\delta}$, we have
  $$
  (1+ \varepsilon) G (z + \zeta)  \leq   \psi_n (z, \mathcal A) + \varepsilon \log (\delta\slash R_0) + C'_0   + C'_1 r \delta^{-1\slash \gamma_0},
  $$
    
 Set  $ \delta  :=r_0^{-\gamma_0}  r^{\gamma_0}$ so that $r \leq r_0 \delta$.  Then we define $ \varepsilon = \varepsilon (r) > 0$ so that 
 $\varepsilon \log (\delta \slash R_0) + C'_0   + C'_1 r^{-\gamma_0}= 0$ i.e. 
 $$
 \varepsilon (r) :=    C_2   \left[\log (R_0\slash r^{\gamma_0})\right]^{-1},
 $$
 where $C_2 :=  C'_0 + C'_1 r^{-\gamma_0}$.
Since $0 < r  \leq r_0 \delta_0^{1\slash \gamma_0} \leq r_0$, we have $\varepsilon (r) \leq C_2   \left[\log (R_0\slash r_0^{\gamma_0})\right]^{-1}$, we choose $\varepsilon_0 :=  C_2   \left[\log (R_0\slash r_0^{\gamma_0})\right]^{-1}$. 

Fix $0 < r \leq r_1$ and $\vert \zeta \vert = r$. Then  for any  $z \in \Omega^\zeta \cap \partial A_{\delta}$,  we have
  $$
  (1+ \varepsilon) G^\zeta (z)   \leq   \psi_n (z, \mathcal A).
  $$
  
  Now we can define as above the function $v^\zeta$ by the formula (\ref{eq:vzeta}) with $m = n$ and the corresponding subextension $w^\zeta$ by (\ref{eq:wzeta}) to obtain the basic  inequality \ref{eq:basic1} with $m=n$ i.e.
  \begin{equation} \label{eq:basic2}
  (1+\varepsilon) G^\zeta  \leq G(z) + 2 M (1+ \varepsilon) r
  \end{equation}
  in $\Omega \cap \Omega^\zeta \setminus A_\delta$.
 
  Therefore   for any   $\vert \zeta \vert \leq r \leq r_1$ and for any $z \in (\Omega^\zeta \cap \Omega) \setminus A_\delta$,   we get 

\begin{eqnarray} \label{eq:Eq1}
G(z+ \zeta)  - G(z) \leq  -  \varepsilon  G(z + \zeta) & =    & - C_2 \left[\log (R_0\slash r^{\gamma_0})\right]^{-1}   G(z + \zeta)
\nonumber \\
&& + 2 (1+\varepsilon) C M  r.
\end{eqnarray}
 
 On the other hand, fix $ \eta > 0 $ and estimate $-  G(z + \zeta) $ from above  when  $z \notin A_{2 \eta}$. Indeed for such $z$ we have $ z' := z + \zeta \notin A_{\eta}$
if  $\vert \zeta \vert = r \leq r_0 \eta$. 
 Then for any $a  \in A$ we have $\vert z+ \zeta  - a\vert \geq \vert z  - a\vert - r \geq \eta$.
 It follows from (\ref{eq:FIneq4}) that
  $$
  -G (z+\zeta) \leq  \sum_{a \in A} - \nu(a) \log (\vert z + \zeta - a\vert \slash R_0)  \leq  \gamma_1 \log(R_0\slash \eta).
   $$
Fix $0 < \alpha< 1$ and  define $\eta > 0$ so that
  $$
   \log (R_0 \slash (r_0\eta)^{\gamma_0}) = [ \log(R_0 \slash r^{\gamma_0})]^{1-\alpha}
  $$
   Then $r \leq r_0 \eta$ and since $r_0 < 1$ and $R_0^{\gamma_0} \leq R_0$, it follows from  the previous inequality that for $\vert \zeta \vert \leq r \leq r_0$ and $z \in \Omega \setminus A_{2 \eta}$, we have
 \begin{equation} \label{eq:Eq2}
  -G (z+\zeta) \leq \gamma_1 \gamma_0^{-1}\log(R_0\slash \eta^{\gamma_0}) .
 \end{equation}
 Since $A_ r \subset A_\eta$, we have  $\Omega \setminus A_\eta \subset \Omega \setminus A_{r}$. Then we can combine the estimates (\ref{eq:Eq1}) and (\ref{eq:Eq2}) to obtain the following estimate:
 for any $z \in \Omega \setminus A_\eta$, $\vert \zeta \vert = r \leq r_2$ and $z \in \Omega^\zeta \setminus A_\eta$, 
 \begin{eqnarray}
G(z+ \zeta) - G(z) &\leq &   C_2  \gamma_1 \gamma_0{-1} [ \log(R_0 \slash r^{\gamma_0})]^{-\alpha}  +  2 (1+\varepsilon_0) C M  r \nonumber \\
  \nonumber \\
& \leq &   C_3 [ \log(R_1 \slash r)]^{-\alpha},
\end{eqnarray}
where $R_1 := R_0^{1 \slash \gamma_0}$ and $C_3 = C_3(R_0,\gamma_0)> 0$ is a uniform constant.

Now we proceed to the proof of the uniform continuity of $\exp G_n (\cdot,\mathcal A, \Omega)$ in $\bar \Omega$ with a uniform control of its modulus of continuity when  $\mathcal A \in \mathcal{F} (\delta_0,\gamma_0,\gamma_1,\sigma_0).$

Observe that for $x, y \in ]-\infty,0]$ we have $\vert e^x - e^y \vert \leq \vert x - y\vert$. Hence  for any $\mathcal{A}  \in \mathcal{F} (\delta_0,\gamma_0,\gamma_1,\sigma_0)$, $z \in \Omega \setminus A_{\eta}$ and $z' \in \Omega \setminus A_\eta$  such that  $\vert z'-z\vert \leq r \leq r_2$, we have 
\begin{equation}\label{eq:Holder-exp1}
 \exp G_n (z', \mathcal A,\Omega) - \exp G_n (z, \mathcal A,\Omega)  \leq C_3 [ \log(R_1 \slash r)]^{-\alpha}.
\end{equation} 

Now we want to estimate the left hand side of the inequality (\ref{eq:Holder-exp1}) for $z \in A_{\eta}$, $z' \in \Omega$  such that  $\vert z'-z\vert \leq r_2$.
Indeed, fix $ \mathcal A  \in \mathcal F (\delta_0,\gamma_0,\gamma_1,\sigma_0)$ and let $z \in A_{\eta}$, $z' \in \Omega$  such that  $\vert z'-z\vert = r \leq r_2$.
Then there exists $(a,\nu) \in \mathcal A$ such that $\vert z - a\vert \leq R_0  (\eta \slash R_0)^{1\slash \nu} \leq R_0  (\eta \slash R_0)^{1\slash \gamma_1}$ and then
 $$ 
 \nu \log(\vert z'-a\vert\slash R_0) \leq \gamma_0 \log [(\eta \slash R_0)^{1\slash \gamma_1}+ r \slash R_0] 
 $$
Then by  (\ref{eq:FIneq4}), we have 
$$
G_n (z',\mathcal A,\Omega) \leq \nu \log( \vert z'-a\vert\slash R_0) + \log (R_0 \slash \delta_0)  \leq \gamma_0 \log [(\eta \slash R_0)^{1\slash \gamma_1}
+ r \slash R_0] + C_0,
$$ 
for any $z\in \Omega$. Hence
$$
 \exp G_n (z', \mathcal A,\Omega)  -  \exp G_n (z, \mathcal A,\Omega)  \leq e^{C_0}  [(\eta \slash R_0)^{1\slash \gamma_1}+ r \slash R_0]^{\gamma_0} =: g(r).
 $$
It is easy to see that  there exists a constant $K = K (m,n,\alpha,\gamma_0,R_0)>0$ such that for $r\leq r_1$, we have
$$ g(r) \leq K [\log R_1\slash r]^{-\alpha}.$$
 Hence for any $z \in A_{\eta}$, $z' \in \Omega$  such that  $\vert z'-z\vert = r \leq r_2$, we have
 \begin{equation}\label{eq:Holder-exp2}
 \exp G_n (z', \mathcal A,\Omega)  -  \exp G_n (z, \mathcal A,\Omega)  \leq L_2  [\log R_1\slash r]^{-\alpha}.
\end{equation} 
 The inequality (\ref{eq:equicontespace}) of the theorem follows from (\ref{eq:Holder-exp1}) and (\ref{eq:Holder-exp2}).
 \end{proof}
 
 \smallskip

\begin{remark} \label{rem:Lipschitz} In the previous result we have assume the domain $\Omega$ to be of Lipschitz type i.e. it admits a bounded Lipschitz continuous exhaustive $m$-subharmonic function. This condition is a strong condition which is different from the condition that $\Omega$ has a Lipschitz boundary. 

 However if we  only assume that the domain $\Omega$ is $m$-hyperconvex, we can still prove that the exponential Green function is uniformly continuous in $\bar \Omega$ and get a precise estimate on its modulus of continuity. 
 Indeed, in this case it is easy to see from the previous proof that the modulus of continuity of $\exp G_m(z, \mathcal A,\Omega)$ will be uniformly controlled by the modulus of continuity defined for $r > 0$ small enough by the following formula :
$$
\omega (r) := \max \{[\log R_1\slash r]^{-\alpha}, \omega_\rho (r)\}, 
$$ 
 where $\omega_\rho (r) $ is the modulus of continuity of the exhaustion function $\rho$ and $\alpha \in ]0,1[$.
 
 This proves in particular the uniform continuity of the exponential Green function when the domain is only assumed to be $m$-hyperconvex which is the result proved by Lelong \cite{Lel89} for the pluricomplex Green function i.e.  in the case $m=n$.
\end{remark}

\subsection{Equicontinuity in the weighted space variable} 
The third step in the proof of Theorem B will consist in proving  the following result.
\begin{theorem} \label{thm:Holder-espacepondere}
 Let $\delta_0 > 0$ be small enough and $ 0 < \gamma_0 < \gamma_1, \sigma_0 > 0$ be fixed constants and  $\mathcal F  = \mathcal F (\delta_0,\gamma_0,\gamma_1,\sigma_0)$. 
 
 Then the family $\{\exp G_m (z; \cdot; \Omega) \, ; \,  z \in \bar{\Omega}\}$ is equicontinuous in the metric space $ (\mathcal F , d_H)$. 
  More precisely we have :
   
   1. If $1 \leq m < n$,  for any $0 < \tau < \frac{2s}{2s+1}$,   for any $0 < \tau < \frac{2 s}{2 s +1}$   there exists constants  $K > 0$ and $r_1 >0$ depending on $ (\tau,\delta_0,\gamma_0,\gamma_1,m,n) $ such that for any $ z \in \bar \Omega $,  and  $ \mathcal A ,  \mathcal A' \in \mathcal F$ with $ d_H (\mathcal A;\mathcal A') \leq r \leq r_1$,  we have
\begin{equation} \label{eq:Holder-espacepondere}
\vert \exp G_m(z, \mathcal A',\Omega) - \exp G(z, \mathcal A,\Omega)\vert \leq  K r^{\tau}.
\end{equation}

2.  If $m=n$, for any $0 < \alpha < 1$,   there exists a constant  $L = L (n,\alpha, \delta_0,\gamma_0,\gamma_1)> 0$ and $r_1 = r_1 (n,\tau,\delta_0,\gamma_0,\gamma_1)> 0$ such that for any $ z \in \bar \Omega $,  and  $ \mathcal A ,  \mathcal A' \in \mathcal F$ with $ d_H (\mathcal A;\mathcal A') \leq r \leq r_1$,  we have 
\begin{equation} \label{eq:equicont-espace}
\vert  \exp G_n(z,\mathcal A', \Omega) - \exp G_n(z,\mathcal A,\Omega) \vert \, \leq \, L \cdot  \left[\log (R_1 \slash r)\right]^{- \alpha},
\end{equation}
where $R_1 := R_0^{1 \slash \gamma_0}$.
\end{theorem}
\begin{proof} As before the proof is divided in two cases. Since $\Omega$ is fixed, we set 
$G_m(z,\mathcal A) := G_m(z,\mathcal A, \Omega)$. 

\smallskip

1. {\it The case $m < n$}.
  Fix  $\mathcal A, \mathcal A' \in \mathcal F (\delta_0,\gamma_0,\gamma_1,\sigma_0)$ and $\varepsilon > 0$. 
Then by (\ref{eq:FIneq2}), (\ref{eq:FIneq3}) and  (\ref{eq:FIneq7}),  for any $0 < \delta\leq \delta_0$ and  $z \in \Omega \setminus (A_{\delta} \cup A'_{\delta})$, we have
  \begin{eqnarray*}
(1+ \varepsilon) G_m(z,\mathcal A') &\leq & (1+ \varepsilon) \phi_m(z,\mathcal A') + (1+ \varepsilon) \delta_0^{-2 s} \\
 &\leq & \psi_m (z, \mathcal A) + \varepsilon \psi_m (z, \mathcal A) +  (1+ \varepsilon)(\delta_0^{-2 s} + \gamma_1^2 \gamma_0^{-1}\sigma_0^{-2s})  \\
 &+ & (1+ \varepsilon) { L'_m} \delta^{-2 s -1} d_H (\mathcal A;\mathcal A').
  \end{eqnarray*}
  Fix $\varepsilon_0 >0$ and set $C_0 := (1+ \varepsilon_0)(\delta_0^{-2 s} + \gamma_1^2 \gamma_0^{-1}\sigma_0^{-2s})$ and $C_1 :=  (1+ \varepsilon_0) L'_m$.
  
  Then for any $0 < \delta\leq \delta_0$, $0 < \varepsilon \leq \varepsilon_0$, and  $z \in \Omega \setminus (A_{\delta} \cup A'_{\delta})$, we have
  $$
  (1+ \varepsilon) G_m(z,\mathcal A') \leq \psi_m (z, \mathcal A) + \varepsilon \psi_m (z, \mathcal A) + C_0 + C_1 \delta^{-2 s -1} d_H (\mathcal A;\mathcal A').
  $$
  
   Recall that 
   $$
   A_\delta = \{z \in \Omega ; \phi_m (z,\mathcal A) < - \delta^{-2 s}\} \, \, \mathrm{and}\, \, A'_\delta = \{z \in \Omega ; \phi_m (z,\mathcal A') < - \delta^{-2 s}\}, 
   $$
 and observe that if  $z \in \partial A_\delta$ then $\phi_m(z,\mathcal A) = - \delta^{-2s}$. On the other hand   if  $z \in \partial A'_\delta$, we have $\psi_m(z,\mathcal A') \leq \phi_m(z,\mathcal A') = -\delta^{-2 s}$ and by Lemma \ref{lem:Flem} it follows that for  $z \in \partial A'_\delta$, we have
 \begin{eqnarray*}
 \psi_m(z,\mathcal A)& \leq & \psi_m(z,\mathcal A')  +   L'_m  \delta^{-2s-1} d_H (\mathcal A,\mathcal A') \\
 &\leq & - \delta^{-2s} +  L'_m  \delta^{-2s-1} d_H (\mathcal A;\mathcal A').
 \end{eqnarray*}
Therefore for any fixed $r >0$, by the previous inequality it follows that  for   $ z \in \partial (A_\delta \cup A'_\delta)$ and $ d_H (\mathcal A;\mathcal A')\leq r$, 
 \begin{eqnarray*}
(1+ \varepsilon) G_m(z,\mathcal A') &\leq& \psi_m (z, \mathcal A) - \varepsilon \delta^{-2 s} + \varepsilon  
L'_m r  \delta^{-2s-1}  + C_0 + C_1 r \delta^{-2 s -1}  \\
&\leq & \psi_m (z, \mathcal A) - \varepsilon \delta^{-2 s} + C_0 + 2 C_1 r \delta^{-2 s -1}.
 \end{eqnarray*}
 Now  set  $\delta = r^{1 \slash (2 s +1)}$ for $0 <  r \leq  r_0:= \delta_0^{1 \slash (2s +1)}$ and choose $\varepsilon >0$ so that 
 $- \varepsilon \delta^{-2 s} + C_0 + 2 C_1 r \delta^{-2 s -1} = 0$ i.e. 
 $$
 \varepsilon = \varepsilon(r) := (C_0   + 2 C_1) \delta^{2 s} = C_2 \delta^{2 s} = C_2 r^{2s \slash (2 s +1)}.
 $$
 If we define $\varepsilon_0 := C_2 r_0^{2s\slash (2s+1)}$. The previous inequality yields  $(1+ \varepsilon) G_m(z,\mathcal A') \leq \psi_m(z,\mathcal A)$ in $\partial (A_\delta \cup A'_\delta)$ for $0< \varepsilon \leq \varepsilon_0$.
 
 Fix $0< r \leq r_0$. Then by the gluing principle, the following function
\begin{equation*} 
 v (z)  :=   \left\{\begin{array}{lcl} 
  \psi_m (z,\mathcal A)  ,  &\hbox{in}\  A_\delta \cup A'_\delta, \\
  \max \{(1 + \varepsilon) G_m(z,\mathcal A') , \psi_m (z, \mathcal A)\}, & \hbox{in}\  \Omega \setminus (A_\delta \cup A'_\delta).
\end{array}\right.
\end{equation*}
is a negative $m$-subharmonic in $\Omega$ such that for any $a\in A$, $\nu_m(v,a) = \nu_m ( \psi_m (\cdot, \mathcal A),a) = \nu(a)$. By  the formula (\ref{eq:LelongFormula}),  it satisfies $v \leq G_m(\cdot,\mathcal A)$ in $\Omega$.  
This implies  that for $z \in \Omega \setminus (A_\delta \cup A'_\delta)$, we have

\begin{equation} \label{eq:G1}
 (1 + \varepsilon) G_m(z,\mathcal A') \leq G_m (z,\mathcal A). 
\end{equation}

Finally for  $0 <  r \leq  r_0$, $\mathcal A, \mathcal A' \in \mathcal F (\delta_0,\gamma_0,\gamma_1,\sigma_0)$ with $d_H(\mathcal A, \mathcal A') \leq r$ and $z \in \Omega \setminus (A_\delta \cup A'_\delta)$, we have 
\begin{equation} \label{eq:Estimate1}
G_m (z,\mathcal A') -  G_m (z,\mathcal A)  \leq - \varepsilon (r) G_m (z,\mathcal A').
\end{equation}

Now fix $0 < \alpha < 1$ and observe that for $z \in \Omega \setminus ( A_{\delta^{\alpha}} \cup  A'_{\delta^{ \alpha}}) \subset  \Omega \setminus  A'_{\delta^{ \alpha}}$,
\begin{eqnarray} \label{eq:Estimate2}
 - G_m (z, \mathcal A') \leq  - \psi_m (z, \mathcal A')  &= & \sum_{(a',\nu')  \in \mathcal A'} \frac{\nu'}{\vert z-a'\vert^{2s}}   \leq    \gamma_1 \gamma_0^{-1} \delta^{- 2 s\alpha}.  
\end{eqnarray}
   Hence, since $A_\delta \subset A_{\delta^{\alpha}} $ for $\delta \leq 1$, we can apply (\ref{eq:Estimate1})  and (\ref{eq:Estimate2}) to deduce that for  any $z \in \Omega \setminus  (A_{\delta^{ \alpha}} \cup  A'_{\delta^{\alpha}})$,  we have
\begin{eqnarray*}
G_m(z, \mathcal A') - G_m(z, \mathcal A)  &\leq & \varepsilon(r) \gamma_1 \gamma_0^{-1}  \delta^{-2 s \alpha}  \\
&\leq & L_2 r^{2s (1-\alpha)\slash (2s+1)},
\end{eqnarray*}
where $L_2:=    \gamma_1 \gamma_0^{-1} C_2$.

Therefore  
for $z \in \Omega \setminus  (A_{\delta^{ \alpha}} \cup  A'_{\delta^{\alpha}})$ and $\mathcal{A} , \mathcal{A'} \in \mathcal{E} (\delta_0,\gamma_0,\gamma_1)$ such that  $d_H (\mathcal A;\mathcal A')  \leq r \leq r_0$,  we have 
\begin{equation} \label{eq:Holder3}
 G_m(z, \mathcal A') - G_m (z, \mathcal A) \leq L_2 r^{\tau(\alpha)},
\end{equation}
where  $ \tau(\alpha) =  {2s (1-\alpha)\slash (2s+1)}$ and $\delta = r^{1\slash (2s+1)}$.

Now we proceed to the proof of Hölder continuity of $\exp G_m (z,\mathcal A)$ in $\bar \Omega \times \mathcal{F}$.

Observe that for $x, y \in ]-\infty,0]$ we have $\vert e^x - e^y \vert \leq \vert x - y\vert$. Then the inequality (\ref{eq:Holder3}) implies that  for any $z \in \Omega \setminus (A_{\delta^\alpha} \cup A'_{\delta^\alpha})$ and any $\mathcal{A} , \mathcal{A'} \in \mathcal{F} $ such that  $d_H (\mathcal A;\mathcal A') \leq  r\leq r_0$, we have 
\begin{equation}\label{eq:Holder-exp}
\exp G_m (z, \mathcal A') - \exp G_m (z, \mathcal A)  \leq L_2 r^{\tau(\alpha)},
\end{equation} 

Now we want to estimate the left hand side of the inequality (\ref{eq:Holder-exp}) for $z \in A_{\delta^{\alpha}} \cup A'_{\delta^{\alpha}}$.

Fix $ \mathcal A,  \mathcal A' \in  \mathcal{F}$ with $d_H (\mathcal A,  \mathcal A') \leq r \leq r_0$.

Assume first that $z \in A'_{\delta^{\alpha}}$. Then there exists $(a',\nu') \in \mathcal A'$ such that $\vert z - a'\vert \leq \nu'^{1\slash 2s} {\delta^{\alpha}}$, which 
by (\ref{eq:FIneq4}) yields $G_m(z,\mathcal A') \leq - \delta^{-2s \alpha} + \delta_0^{-2s} = -  r^{-2s \alpha \slash (2s+1)} + \delta_0^{-2s}$.
 Hence
$$
 \exp G_m(z, \mathcal A')  -  \exp G_m(z, \mathcal A)  \leq \exp \left(-  r^{-2s \alpha \slash (2s+1)} + \delta_0^{-2s}\right).
$$
Now assume that $z \in A_{\delta^{\alpha}}$. Then  there exists $(a,\nu) \in \mathcal A$ such that  $\vert z - a\vert \leq \nu^{1 \slash 2s} \delta^{\alpha}$.  Since $d_H (\mathcal A,  \mathcal A') \leq r$ there exists $(a',\nu') \in \mathcal A'$ such that 
$\vert a - a'\vert + \vert \nu -  \nu' \vert \leq r$. Moreover $\vert z-a' \vert \leq \vert z - a\vert + \vert a- a' \vert \leq  \gamma_1^{1\slash 2s} {\delta^{\alpha}} + r$.  Hence $G_m (z, \mathcal A') \leq  - \gamma_0 ( \gamma_1^{1\slash 2s} {\delta^{\alpha}} + r)^{-2s} + \delta_0^{-2s}$.
This implies that  
$$
\exp G_m (z, \mathcal A') - \exp G_m (z, \mathcal A) \leq \exp \left(- \gamma_0 ( \gamma_1^{1\slash2s} r^{\alpha\slash(2s +1)} + r)^{-2s} + \delta_0^{-2s}\right).
$$
Finally  define the following modulus of continuity
\begin{eqnarray*}
h (r) :=  \exp \left(-  r^{-2s \alpha \slash (2s+1)}   + \delta_0^{-2s}\right) +  \exp \left(- \gamma_0 ( \gamma_1^{1\slash 2s} r^{\alpha\slash(2s +1)}  + r)^{-2s} + \delta_0^{-2s}\right),
\end{eqnarray*}
and observe  that for any  $z \in A_{\delta^{\alpha}} \cup A'_{\delta^{\alpha}}$  and any $ \mathcal A,  \mathcal A' \in  \mathcal{F} $ with $d_H (\mathcal A,  \mathcal A') \leq r$, we have
\begin{eqnarray} \label{eq:ExpCont}
\exp G_m (z, \mathcal A') - \exp G_m (z, \mathcal A)&\leq &  h(r).
\end{eqnarray}

It's easy to see that there exists   a uniform constant $K_0 = K_0(m,n,\alpha)  > 0$ such that $h (r) \leq K_0 \, r^{\tau(\alpha)}$,  for $0 < r \leq r_0$ and then it follows from (\ref{eq:Holder-exp}) and (\ref{eq:ExpCont}) that for any $z \in \bar{\Omega}$ and   any $ \mathcal A,  \mathcal A' \in  \mathcal{F} $ with $d_H (\mathcal A;\mathcal A') \leq r_0$ we have
\begin{equation} 
\exp G_m (z, \mathcal A') - \exp G_m (z, \mathcal A) \leq  K d_H (\mathcal A;\mathcal A')^{\tau(\alpha)},
\end{equation}
where $K := L_2 + K_0$.
Since $0 < \alpha < 1$ can be taken arbitrarily close to $0$, $\tau (\alpha)  = 2s (1-\alpha)\slash (2s+1)$ can be taken arbitrarily close to $2s \slash (2s+1)$, this implies the inequality (\ref{eq:Holder-espacepondere}). 

\smallskip
 
2. {\it The case $m = n$}.
  Fix  $\mathcal A, \mathcal A' \in \mathcal F (\delta_0,\gamma_0,\gamma_1,\sigma_0)$ and $\varepsilon > 0$. 
Then by (\ref{eq:FIneq2}), (\ref{eq:FIneq3}) and  (\ref{eq:FIneq7}),  for any $0 < \delta\leq \delta_0$ and  $z \in \Omega \setminus (A_{\delta} \cup A'_{\delta})$, we have
 \begin{eqnarray*}
(1+ \varepsilon) G_n(z,\mathcal A') &\leq & (1+ \varepsilon) \phi_n(z,\mathcal A') + (1+ \varepsilon) (\log(R_0\slash\delta_0)  \\
 &\leq & \psi_n (z, \mathcal A) + \varepsilon \psi_n (z, \mathcal A) +  (1+ \varepsilon)(\log(R_0\slash\delta_0) \\
 &+& \gamma_1^2 \gamma_0^{-1} \log(R_0\slash\delta_0) +  (1+ \varepsilon) { L'_n} \delta^{-1 \slash\gamma_0} d_H (\mathcal A;\mathcal A').
  \end{eqnarray*}
  Fix $\varepsilon_0 >0$ and set $C_0 := (1+ \varepsilon_0)(\log(R_0\slash\delta_0) + \gamma_1^2 \gamma_0^{-1} \log(R_0\slash\delta_0))$ and $C_1 :=  (1+ \varepsilon_0) L'_n$.
  
  Then for any $0 < \delta\leq \delta_0$, $0 < \varepsilon \leq \varepsilon_0$, and  $z \in \Omega \setminus (A_{\delta} \cup A'_{\delta})$, we have
  $$
  (1+ \varepsilon) G_n(z,\mathcal A') \leq \psi_n (z, \mathcal A) + \varepsilon \psi_n (z, \mathcal A) + C_0 + C_1 \delta^{-1 \slash\gamma_0}  d_H (\mathcal A;\mathcal A').
  $$
  
   Recall that 
   $$
   A_\delta = \{z \in \Omega ; \phi_n (z,\mathcal A) < \log (\delta \slash R_0)\}, 
   $$
 and observe that if  $z \in \partial A_\delta$ then $\psi_n (z,\mathcal A) \leq \phi_n (z,\mathcal A) = \log (\delta \slash R_0)$. 
 
 Therefore it follows from the previous inequality that  for any fixed $r >0$, and    $ z \in \partial (A_\delta \cup A'_\delta)$ and $ d_H (\mathcal A;\mathcal A')\leq r$, 
 \begin{eqnarray*}
(1+ \varepsilon) G_n (z,\mathcal A') 
&\leq & \psi_n (z, \mathcal A) - \varepsilon \log(R_0\slash \delta) + C_0 +  C_1 r  \delta^{-1 \slash\gamma_0} .
 \end{eqnarray*}
  Now  set  $\delta = r^{\gamma_0}$ for $0 <  r \leq  r_0:= \delta_0^{1 \slash \gamma_0}$ and choose $\varepsilon >0$ so that 
 $$
 - \varepsilon \log(R_0\slash \delta) + C_0 +  C_1 = 0 \, \, \hbox{i.e.} \, \, \, 
 \varepsilon = \varepsilon(r) :=  C_2  \slash \log(R_1 \slash r),
 $$
 where $C_2 := \gamma_0^{-1} (C_0 +C_1)$ and $R_1 := R_0^{1 \slash \gamma_0}$.

 If we define $\varepsilon_0 := C_2  \slash \log(R_1 \slash r_0) $. The previous inequality yields  $(1+ \varepsilon) G_n (z,\mathcal A') \leq \psi_n (z,\mathcal A)$ in $\partial (A_\delta \cup A'_\delta)$ for $0< \varepsilon \leq \varepsilon_0$.
 
 The same gluing process as in the previous case yields for  $0 <  r \leq  r_0$, $\mathcal A, \mathcal A' \in \mathcal F (\delta_0,\gamma_0,\gamma_1,\sigma_0)$ with $d_H(\mathcal A, \mathcal A') \leq r$ and $z \in \Omega \setminus (A_\delta \cup A'_\delta)$, 
\begin{equation} \label{eq:Est1}
G_n (z,\mathcal A') -  G_n (z,\mathcal A)  \leq - \varepsilon (r) G_m (z,\mathcal A').
\end{equation}

Now fix $\delta < \eta < \delta_0 $ and observe that for $z \in \Omega \setminus ( A_{\eta} \cup  A'_{\eta}) \subset  \Omega \setminus  A'_{\delta}$,
\begin{eqnarray} \label{eq:Est2}
 - G_n (z, \mathcal A') \leq  - \psi_n (z, \mathcal A')  &= & - \sum_{(a',\nu')  \in \mathcal A'}  \nu' \log (\vert z-a' \vert  \slash R_0) \nonumber \\
 & \leq &   \gamma_1 \gamma_0^{-1}\log(R_0\slash \eta).  
\end{eqnarray}
   Hence, since $A_\delta \subset A_{\eta} $ for $\delta \leq 1$, we can apply (\ref{eq:Est1})  and (\ref{eq:Est2}) to deduce that for  any $z \in \Omega \setminus  (A_{\eta} \cup  A'_{\eta})$,  we have
\begin{eqnarray*}
G_n (z, \mathcal A') - G_n (z, \mathcal A)  &\leq & \varepsilon(r) \gamma_1 \gamma_0^{-1} \log(R_0\slash \eta) \\
&\leq & L_2 \frac{\log(R_0\slash \eta)}{ \log(R_0 \slash r^{\gamma_0})},
\end{eqnarray*}
where $L_2:=    \gamma_1 \gamma_0^{-1} C_2$.

Now fix $0 < \alpha< 1$ and choose $\eta$ so that
$$
\log(R_0\slash \eta) = \left[\log(R_0 \slash r^{\gamma_0})\right]^{(1 - \alpha)}. 
$$

Then $\delta = r^{\gamma_0} \leq \eta$ and thus  
for $z \in \Omega \setminus  (A_{\eta} \cup  A'_{\eta})$ and $\mathcal{A} , \mathcal{A'} \in \mathcal{E} (\delta_0,\gamma_0,\gamma_1)$ such that  $d_H (\mathcal A;\mathcal A')  \leq r \leq r_0$,  we have 
\begin{equation} \label{eq:Holder3}
 G_n(z, \mathcal A') - G_n(z, \mathcal A) \leq \frac{L'_2}{\left[\log(R_1 \slash r)\right]^{\alpha}},
\end{equation}
where  $L'_2 := \gamma_0^{-\alpha} L_2 $.

To prove the same inequality for $\exp G_n (z,\mathcal A)$ in $\bar \Omega \times \mathcal{F}$ we proceed exactly as in the previous case (see also the proof of Theorem \ref{thm:Holder-espace}).
 \end{proof}
 
 \smallskip
 
 \subsection{Proof of Theorem B}
 Now we can easily deduce Theorem B  from Theorem \ref{thm:Holder-espace} and Theorem \ref{thm:Holder-espacepondere}.
  Indeed, for any $(z,\mathcal A) \in  \Omega \times \mathcal F (\delta_0,\gamma_0,\gamma_1,\sigma_0) $ and $(z',\mathcal A') \in  \Omega \times \mathcal F (\delta_0,\gamma_0,\gamma_1,\sigma_0)$, we have  
 \begin{eqnarray*}
 \vert \exp G_m(z', \mathcal A') - \exp G_m (z, \mathcal A)\vert & \leq  & \vert \exp G_m(z', \mathcal A') - \exp G_m (z, \mathcal A')\vert \\
 & +&   \vert \exp G_m (z, \mathcal A') - \exp G_m (z, \mathcal A)\vert.
 \end{eqnarray*}
 Then the required inequality (\ref{eq:Hldercontinuite}) follows from the inequalities (\ref{eq:Holder-espace}) and (\ref{eq:Holder-espacepondere}).
 
 \begin{remark} The previous result shows that for any fixed $\mathcal A \subset \Omega \times ]0,+\infty[$, the weighted Green function $G_m(\cdot,\mathcal A, \Omega)$ is locally Lipschitz in $\Omega \setminus A$ when $m < n$ and continuous in $\Omega \setminus A$ when $m=n$. More precisely  for a fixed open subset $D \Subset \Omega$, let us denote by  $\mathcal F_D = \mathcal F_D(\delta_0,\gamma_0,\gamma_1,\sigma_0)$  the family of weighted sets $\mathcal A \in \mathcal F (\delta_0,\gamma_0,\gamma_1,\sigma_0)$ such that $\mathcal A \subset D \times ]0,+\infty[$ endowed with the Hausdorff distance. Then the weighted Green function $G_m(\cdot, \cdot, \Omega)$ is  Lipschitz continuous in $(\bar \Omega \setminus D) \times \mathcal F_D$ when $m < n$ and uniformly continuous in $(\bar \Omega \setminus D) \times \mathcal F_D$ when $m=n$.
\end{remark}

\begin{remark} \label{rem:orderedpoles} Theorem B generalizes and improves  the result of Lelong in many directions.
First Lelong considered ordered sets of weighted poles with the same cardinality $p \geq 1$ and proved the uniform continuity for the euclidean distance on the weighted poles. More precisely,  if $\mathcal A := \{(a_k,\nu_k)_{1 \leq k \leq p }\} \subset (\Omega \times ]0,+  \infty[)^p$ and $\mathcal A' := \{(a'_k,\nu'_k)_{1 \leq k \leq p} \} \subset (\Omega \times ]0,+  \infty[)^p$, the distance considered by Lelong is defined as follows :
$$
d_L (\mathcal A, \mathcal A') :=  \sum_{1 \leq k \leq p} (\vert a_k - a'_k\vert + \vert \nu_k - \nu'_k\vert).
$$
This distance is sensitive to the order of the poles, while the Hausdorff distance defined by the formula  (\ref{eq:Hdistance}) is not. 

Moreover, it follows from the definitions  that for any $\mathcal A , \mathcal A'  \subset (\Omega \times ]0,+  \infty[)^p$, we have
$$
d_H (\mathcal A,\mathcal A') \leq d_L (\mathcal A,\mathcal A').
$$
However the two distances are not equivalent. Indeed, according to Lelong's notations, we define the following sequence of sets
$$\mathcal A_j := \{(a_1^j,1), (a_2^j,1)\},   a_1^j := (0,2^{-j}) \in \B^2,  a_2^j := (1 \slash 2,2^{-j}) \in \B^2,
$$
 for $j \in \N^*$.
 
We also define 
$$
\mathcal A := \{(a_1,1), (a_2,1)\}, a_1 := (1\slash 2,0), a_2 := (0,0).
$$
Then it's easy to see that for any $j \in \N$, we have $d_L (A_j,A) = 1+ 2^{-j + 1}  \geq 1$ while $d_H(\mathcal A_j,\mathcal A) = 2^{-j} \to 0$ as $j \to + \infty$.
\end{remark}

\smallskip

 {\bf Dedication :}  This article is dedicated to the memory of Urban Cegrell who made profound contributions in Pluripotential Theory, generalizing the celebrated theory of Bedford and Taylor. His work and especially his famous article \cite{Ceg98} has been  a great source of inspiration for us. 

 \bigskip

{\bf Acknowledgements.} This work was carried out during various visits by the first author to the Institut de Mathématiques de Toulouse (IMT) in spring 2019 and autumn 2021. She would like to thank the IMT for having welcomed her and offering her excellent research conditions. She also thanks the University of Sousse for the financial support that made this collaboration possible.

The authors warmly thank the two reviewers for their careful reading of the first version of this paper. Their insightful comments and suggestions helped to clarify some proofs and to improve the presentation of this version.

\end{document}